\RequirePackage{snapshot}
\makeatletter
\def\snap@providesfile#1[#2]{%
  \wlog{File: #1 #2}%
  \if\expandafter\snap@graphic@test\expanded{#2}@@\@nil
    \snap@record@graphic#1\relax #2 (type ??)\@nil
  \else
    \expandafter\xdef\csname ver@#1\endcsname{#2}%
  \fi
  \endgroup
}
\makeatother

\documentclass[twoside]{book}

\usepackage{lectures}
\usepackage[colorlinks=true,
citecolor=black,
linkcolor=black,
anchorcolor=black,
filecolor=black,
menucolor=black,
urlcolor=black,
pdftitle={Metric geometry on manifolds: two lectures},
pdfsubject={Geometry},
pdfauthor={Anton Petrunin}
]{hyperref}
\makeindex

\begin{document}

\title{Metric geometry on manifolds:
\\ two lectures}
\author{Anton Petrunin}
\date{}
\maketitle

We discuss Besicovitch inequality, width, and systole of manifolds.
I assume that students are familiar with 
measure theory,
smooth manifolds,
degree of map, 
CW-complexes and related notions.

This is a polished version of two final lectures of a graduate course given at Penn State, Spring 2020.
The complete lectures can be found on the author's website;
it includes an introduction to metric geometry \cite{petrunin2020pure}
and elements of Alexandrov geometry based on \cite{alexander-kapovitch-petrunin-2019}.

\parbf{Acknowledgments.}
I want to thank Alexander Lytchak and Alexander Nabutovsky for help.

\thispagestyle{empty}
\tableofcontents
\thispagestyle{empty}

%%%%%%%%%%%%%%%%%%%%%%%%%%%%
%\include{hwa}
%%%%%%%%%%%%%%%%%%%%%%%%%%%%

\chapter{Besicovitch inequality} 

We will focus on Riemannian spaces --- these are specially nice length metrics on manifolds.
These spaces are also most important in applications.

As it will be indicated in Section~\ref{sec:hausdorff-measure},
most of the statements of this and the following lecture have counterparts for general length metrics on manifolds.

\section{Riemannian spaces}

Let $M$ be a smooth connected manifold.
A \index{metric tensor}\emph{metric tensor} on $M$ is a choice of positive definite quadratic forms $g_p$ on each tangent space $\T_pM$ that depends continuously on the point;
that is, in any local coordinates of $M$ the components of $g$ are continuous functions.

A \index{Riemannian!manifold}\emph{Riemannian manifold} $(M,g)$ is a smooth manifold $M$ with a choice of metric tensor $g$ on it.

The \index{length}\emph{$g$-length} of a Lipschitz curve $\gamma\:[a,b]\to M$  is defined by
\[\length_g\gamma=\int_a^b\sqrt{g(\gamma'(t),\gamma'(t))}\cdot dt.\]
The $g$-length induces a metric metric on $M$; it is defined as the greatest lower bound to lengths of Lipschitz curves connecting two given points;
the distance between a pair of points $x,y\in M$ will be denoted by 
\[\dist{x}{y}{g}\quad\text{or}\quad\distfun_x(y)_g.\]
The corresponding metric space $\spc{M}$ will be called \index{Riemannian!space}\emph{Riemannian space}.

\begin{thm}{Exercise}\label{ex:non-differentiable}
Show that isometry between Riemannian spaces might be not induced by a diffeomorphism.

Moreover, there is a continuous Riemannian metric $g$ on $\RR^2$ such that the corresponding Riemannian space admits an isometry to the Euclidean palne but the induced map $\iota\:\RR^2\to\RR^2$ is not differentiable at some point.
\end{thm}

The exercise above shows that in general the smooth structure is not uniquely defined on Riemannian space.
Therefore in general case one has to distinguish between Riemannian manifold and the corresponding Riemannian space altho there is almost no difference.%
\footnote{In fact a straightforward smoothing procedure shows that isometry between Riemannian spaces can be approximated by diffeomorphisms between underlying manifolds; in particular these manifolds are diffeomorphic.
Also, if the metric tensor is smooth, then it is not hard to show that Riemannian space {}\emph{remembers} everything about the Riemannian manifold, in particular the smooth structure;
it is a part of the so-called Myers--Steenrod theorem \cite{myers-steenrod}.} 

The following observation states the key property of Riemannian spaces;
it will be used to extend results from Euclidean space to Riemannian spaces.

\begin{thm}{Observation}\label{obs:lip-chart}
For any point $p$ in a Riemannian space $\spc{M}$ and any $\eps>0$ there is a $e^{\mp\eps}$-bilipschitz chart $s\:W\to V$ from an open subset $W$ of the $n$-dimensional Euclidean space to some neighborhood $V\ni p$.
\end{thm}

\parit{Proof.}
Choose a chart $s\:U\to \spc{M}$ that covers $p$.
Note that there is a linear transformation $L$ such that for the metric tensor in the chart $s\circ L$ is coincides with the standard Euclidean tensor at the point $x=(s\circ L)^{-1}(p)$.

Since the metric tensor is continuous, the restriction of $s\circ L$ to a small neighborhood of $x$ is $e^{\mp\eps}$-bilipschitz.
\qeds

\section{Volume and Hausdorff measure}\label{sec:vol-haus}

Let $(M,g)$ be an $n$-dimensional Riemannian manifold.
If a Borel set $R\subset M$ is covered by one chart $\iota\:U\to M$,
then its \index{volume}\emph{volume} (briefly, $\vol R$ or $\vol_n R$) is defined by 
\[\vol R
\df
\int_{\iota^{-1}(R)}\sqrt{\det{g}}.\]
In the general case we can subdivide $R$ into a countable collection of regions $R_1,R_2\dots$ such that each region $R_i$ is covered by one chart $\iota_i\:U_i\to M$ and define
\[\vol R\df \vol R_1+\vol R_2+\dots\]
The chain rule for multiple integrals implies that the right-hand side does not depend on the choice of subdivision and the choice of charts.

Similarly, we define integral along $(M,g)$.
Any Borel function $u\:M\z\to \RR$, can be presented as a sum $u_1+u_2+\cdots$ such that the support of each function $u_i$ can be covered by one chart $\iota_i\:U_i\to M$
and set 
\[\int_{p\in\spc{M}} u(p)
\df
\sum_i\left[\int_{x\in U_i} u_i\circ s(x)\cdot\sqrt{\det{g}}\right].
\]
In particular
\[\vol R=\int_{p\in R} 1.\]

Let $\spc{X}$ be a metric space and $R\subset \spc{X}$.
The \index{Hausdorff measure}\emph{$\alpha$-dimensional Hausdorff measure} of $R$ is defined by 
$$\haus_\alpha R
\df
\lim_{\eps\to0}
\,
\inf
\set{\sum_{n\in\NN}(\diam A_n)^\alpha}
{\begin{aligned}
&\diam A_n<\eps\ \text{for}
\\
&\text{for each}\ n,\text{all}\  A_n
\\
&\text{are closed, and} 
\\
& \bigcup_{n\in\NN}A_n\supset R.
\end{aligned}
}.$$
For properties of Hausdorff measure we refer to the classical book of  Herbert Federer \cite{federer};
in particular, $\haus_\alpha$ is indeed a measure and $\haus_\alpha$-measurable sets include all Borel sets.

The following observation follows from \ref{obs:lip-chart} and Rademacher's theorem:

\begin{thm}{Observation}\label{obs:lipcart+}
Suppose that a Borel set $R$ in an $n$-dimensional Riemannian space $\spc{M}$ is subdivided into a countable collection of subsets $R_i$ such that each $R_i$ is covered by an $e^{\mp\eps}$-bilipschitz charts
$s_i$.
Then
\begin{align*}
\vol_n R&\lege e^{\pm n\cdot\eps}\cdot\sum_i\vol_n[s_i^{-1}(R_i)]
\intertext{and}
\haus_n R&\lege e^{\pm n\cdot\eps}\cdot\sum_i\haus_n[s_i^{-1}(R_i)]
\end{align*}

\end{thm}

According to \index{Haar's theorem}\emph{Haar's theorem}, 
a measure on $n$-dimensional Euclidean space that is invariant with respect to parallel translations is proportional to volume.
Observe that 
\begin{itemize}
\item A ball in $n$-dimensional Euclidean space of diameter $1$ has unit Hausdorff measure.
\item A unit cube in $n$-dimensional Euclidean space has unit volume.
\end{itemize}
Therefore, for any Borel region $R\subset \EE^n$, we have 
\[\vol_n R=\tfrac{\omega_n}{2^n}\cdot\haus_n R,\eqlbl{eq:vol/mu}\]
where $\omega_n$ denotes the volume of a unit ball in the $n$-dimensional Euclidean space.

Applying \ref{eq:vol/mu} together with \ref{obs:lipcart+}, we get that the inequalities
\[\vol_n R\lege e^{\pm2\cdot n\cdot \eps}\cdot\tfrac{\omega_n}{2^n}\cdot\haus_n R\]
hold for any $\eps>0$.
Since $\eps>0$ is arbitrary, we get that \ref{eq:vol/mu} holds in $n$-dimensional Riemannian spaces.
More precisely:

\begin{thm}{Proposition}\label{prop:vol=haus}
The identity 
\[\vol_n R=\tfrac{\omega_n}{2^n}\cdot\haus_n R\]
holds for any Borel region $R$ in an $n$-dimensional Riemannian space. 
\end{thm}

Since the Hausdorff measure is defined in pure metric terms, the proposition gives another way to prove that the volume does not depend on the choice of chars and subdivision of $R$.

The identity in this proposition will be used to define \index{volume}\emph{volume of any dimension}.
Namely, given an integer $k\ge 0$, the $k$-volume is defined by
\[\vol_k\df\tfrac{\omega_k}{2^k}\cdot\haus_k.\]
By \ref{prop:vol=haus}, if $A$ is a subset of $k$-dimensional submanifold $\spc{N}\subset \spc{M}$, then the two definitions of $\vol_kA$ agree; but the latter definition works for a wider class of sets. 

\begin{thm}{Exercise}\label{ex:volume-preserving+short}
Let $f\:\spc{M}\to \spc{N}$ be a short volume-preserving map between $n$-dimensional Riemannian spaces.
Prove the following statements and use them to conclude that $f$ is locally distance-preserving.

\begin{subthm}{ex:volume-preserving+short:injective}
$f$ is injective; 
that is, if $f(x)=f(y)$, then $x=y$.
\end{subthm}

\begin{subthm}{ex:volume-preserving+short:bi}
For any $c<1$, the map $f$ is locally $[c,1]$-bilipschitz;
that is, for any point in $\spc{M}$ there is a neighborhood $\Omega$ and $\eps>0$ such that the inequality 
\[c\le \frac{|f(x)-f(y)|_{\spc{N}}}{|x-y|_{\spc{M}}}\le 1 \]
holds for any pair of distinct points $x,y\in \Omega$.
\end{subthm}

\end{thm}

\section{Area and coarea formulas}

Suppose that $f\:\spc{M}\to\spc{N}$ is a Lipschitz map between $n$-dimensional Riemannian spaces $\spc{M}$ and $\spc{N}$.
Then by \index{Rademacher's theorem}\emph{Rademacher's theorem} 
the differential $d_p f\:\T_p\spc{M}\to\T_{f(p)}\spc{N}$ is defined at \index{almost all}\emph{almost all} $p\in \spc{M}$;
that is, the differential defined at all points $p\in\spc{M}$ with exception of a subset with vanishing volume.

The differential is a linear map; it defines the Jacobian matrix $\Jac_pf$ in orthonormal frames of $\T_p$ and $\T_{f(p)}\spc{N}$.
The determinant of $\Jac_pf$ will be denoted by $\jac_p$.
Note that the absolute value $|\jac_p|$ does not depend on the choice of the orthonormal frames.

The identity in the following proposition is called \index{area formula}\emph{area formula}.

\begin{thm}{Proposition}
Let $f\:\spc{M}\to\spc{N}$ be a Lipschitz map between $n$-dimensional Riemannian spaces $\spc{M}$.
Then for  any Borel function $u\:\spc{M}\z\to \RR$ the following equality holds:
\[\int_{p\in \spc{M}} u(p)\cdot |\jac_pf|=\int_{q\in \spc{N}}\sum_{p\in f^{-1}(q)} u(p).\]

\end{thm}

\parit{Proof.}
If $\spc{M}$ and $\spc{N}$ are isometric to the $n$-dimensional Euclidean space, then the statement follows from the standard area formula \cite[3.2.3]{federer}.

Note that Jacobian of a $e^{\mp\eps}$-bilipschitz map between $n$-dimensional Riemannian manifolds (if defined) has determinant in the range $e^{\mp n\cdot\eps}$.
Applying \ref{obs:lipcart+} and the area formula in $\EE^n$, we get the following approximate version of the needed identity for any $u\ge0$: 
\[\int_{p\in \spc{M}} u(p)\cdot |\jac_pf|
\lege e^{\pm 3\cdot n\cdot \eps}\int_{q\in \spc{N}}\sum_{p\in f^{-1}(q)} u(p).\]

Since $\eps>0$ is arbitrary, we get that the area formula holds if $u\ge 0$.
Finally, since both sides of the area formula are linear in $u$, it holds for any $u$.
\qeds

The following inequality is called \index{area inequality}\emph{area inequality}:

\begin{thm}{Corollary}\label{cor:area-inequality}
Let $f\:\spc{M}\to\spc{N}$ be a locally Lipschitz map between $n$-dimensional Riemannian spaces.
Then 
\[\int_{p\in A} |\jac_p f|\ge \vol[f(A)]\]
for any Borel subset $A\subset M$.

In particular, if $|\jac_p f|\le 1$ almost everywhere in $A$, then 
\[\vol A \ge \vol[f(A)].\]
\end{thm}

\parit{Proof.} Apply the area formula to the characteristic function of $A$.
\qeds

Suppose that $f\:\spc{M}\to\RR$ is a Lipschitz function defined on an $n$-dimensional Riemannian space $\spc{M}$.
Then by Rademacher's theorem, the differential $d_pf\:\T_p\spc{M}\to\RR$  and the gradient 
$\nabla_pf\in\T_p\spc{M}$ are defined at almost all $p\in \spc{M}$.

The identity in the following proposition is a partial case of the so-called \index{coarea formula}\emph{coarea formula}.
(The general coarea formula deals with the maps to the spaces of arbitrary dimension, not necessary $1$.)

\begin{thm}{Proposition}\label{prop:coarea}
Let $f\:\spc{M}\to\RR$ be a Lipschitz function defined on an $n$-dimensional Riemannian space $\spc{M}$.
Suppose that the level sets $L_x\df f^{-1}(x)$ are equipped with $(n-1)$-dimensional volume $\vol_{n-1}\z\df\tfrac{\omega_{n-1}}{2^{n-1}}\cdot \haus_{n-1}$.
Then for any Borel function $u\:\spc{M}\to \RR$ the following equality holds
\[\int_{p\in \spc{M}} u(p)\cdot |\nabla_pf|=\int_{-\infty}^{+\infty} \left(\,\int_{p\in L_x} u(p)\,\right)\cdot dx.\]
\end{thm}

The following corollary is a partial case of the so-called  \index{coarea inequality}\emph{coarea inequality};

\begin{thm}{Corollary}\label{cor:coarea}
Let $\spc{M}$, $f$, and $L_x$ be as in \ref{prop:coarea}.

Suppose that $f$ is 1-Lipschitz.
Then for any Borel subset $A\subset M$ we have
\[\vol_n A\ge \int_{x\in\RR} \vol_{n-1}[A\cap L_x]\cdot dx.\eqlbl{eq:coarea-inq}\]
\end{thm}

The right-hand side in \ref{eq:coarea-inq} is called \index{coarea}\emph{coarea of the restriction $f|_A$}.

\parit{Instead of proof of \ref{prop:coarea} and \ref{cor:coarea}.}
If $\spc{M}$ is isometric to Euclidean space, then the statement follows from the standard coarea formula \cite[3.2.12]{federer}.
The reduction to the Euclidean space is done the same way as in the proof of the area formula.

To prove the corollary, choose $u$ to be the characteristic function of $A$ and apply the coarea formula.
\qeds

\section{Besicovitch inequality}

A closed connected region in a Riemannian manifold bounded by hypersurface will be called \index{Riemannian!manifold with boundary}\emph{Riemannian manifold with boundary}.
We always assume that the hypersurface can be realized locally as a graph of Lipschitz function in a suitable chart.
In this case one can define $g$-length, $g$-distance, and $g$-volume the same way as we did for usual Riemannian manifolds.

\begin{thm}{Exercise}\label{ex:compact-interior}
Suppose that $(M,g)$ is a compact Riemannian manifold with boundary. 
Observe that the interior $(M^\circ,g)$ of $(M,g)$ is a usual Riemannian manifold.
Show that the space of $(M,g)$ is isometric to the completion of the space of $(M^\circ,g)$.
\end{thm}

\begin{thm}{Theorem}\label{thm:besikovitch}
Let $g$ be a continuous metric tensor on a unit $n$-dimensional cube $\square$.
Suppose that the $g$-distances between the opposite faces of $\square$ are at least $1$; that is, any Lipschitz curve that connects opposite faces has $g$-length at least $1$.
Then \[\vol(\square, g)\ge 1.\]

\end{thm}

This is a partial case of the theorem proved by Abram Besicovitch \cite{besicovitch}.

\parit{Proof.}
We will consider the case $n=2$; the other cases are proved the same way.

\begin{wrapfigure}{r}{30mm}
\vskip-0mm
\centering
\includegraphics{mppics/pic-1320}
\end{wrapfigure}

Denote by $A$, $A'$, and $B$, $B'$ the opposite faces of the square~$\square$.
Consider two functions
\begin{align*}
f_A(x)&\df\min\{\,\distfun_A(x)_g,1\,\},
\\
f_B(x)&\df\min\{\,\distfun_B(x)_g,1\,\}.
\end{align*}
Let $\bm{f}\:\square\to\square$ be the map with coordinate functions $f_A$ and $f_B$;
that is, $\bm{f}(x)\df(f_A(x), f_B(x))$.

\begin{clm}{}\label{f:A->A}
The map $\bm{f}$ sends each face of $\square$ to itself.
\end{clm}

Indeed, 
\[x\in A \quad\Longrightarrow\quad \distfun_A(x)_g=0 \quad\Longrightarrow\quad f_A(x)=0 \quad\Longrightarrow\quad \bm{f}(x)\in A.\]
Similarly, if $x\in B$, then $\bm{f}(x)\in B$.
Further, 
\[x\in A'
\quad\Longrightarrow\quad 
\distfun_A(x)_g\ge 1 
\quad\Longrightarrow\quad 
f_A(x)=1 
\quad\Longrightarrow\quad 
\bm{f}(x)\in A'.\]
Similarly, if $x\in B'$, then $\bm{f}(x)\in B'$.

By \ref{f:A->A}, it follows 
\[\bm{f}_t(x)= t\cdot x + (1-t)\cdot \bm{f}(x)\]
defines a homotopy of maps of the pair of spaces $(\square,\partial \square)$ from $\bm{f}$ to the identity map;
that is, $(t,x)\mapsto \bm{f}_t(x)$ is a continuous map and if $x\in \partial \square$, then $\bm{f}_t(x)\in \partial \square$ for any $t\in [0,1]$.

It follows that $\deg\bm{f}=1$; that is, $\bm{f}$ sends the fundamental class of $(\square,\partial \square)$ to itself.%
\footnote{Here and further, we assume that homologies are taken with the coefficients in $\ZZ_2$, but you are welcome to play with other coefficients.}
In particular $\bm{f}$ is onto.

Suppose that Jacobian  matrix $\Jac_p\bm{f}$ of $\bm{f}$ is defined at $p\in \square$.
Choose an orthonormal frame in $\T_p$ with respect to $g$ and the standard frame in the target $\square$.
Observe that the differentials $d_pf_A$ and $d_pf_B$ written in these frames are the rows of $\Jac_p\bm{f}$.
Evidently $|d_pf_A|\le 1$ and $|d_pf_B|\le 1$.
Since the determinant of a matrix is the volume of the parallelepiped spanned on its rows, we get 
\[|\jac_p \bm{f}|\le |d_pf_A|\cdot|d_pf_B|\le 1.\]
Since $\bm{f}\:\square\to\square$ is a Lipschitz onto map, the area inequality (\ref{cor:area-inequality}) implies that 
\[\vol(\square,g)\ge \vol\square=1.\]
\qedsf

If the $g$-distances between the opposite sides are $d_1,\dots ,d_n$, then following the same lines  one get that 
$\vol (\square,g)\ge d_1\cdots d_n$.
Also note that in the proof we use topology of the $n$-cube only once, to show that the map $f$ has degree one.
Taking all this into account we get the following generalization of \ref{thm:besikovitch}:

\begin{thm}{Theorem}\label{thm:besikovitch+}
Let $(M,g)$ be an $n$-dimensional Riemannian manifold with coonected boundary $\partial M$.
Suppose that there is a degree 1 map $\partial M\to \partial\square$;
denote by $d_1,\dots, d_n$ the $g$-distances between the inverse images of pairs of opposite faces of $\square$ in $M$.
Then 
\[\vol(M,g)\ge d_1\cdots d_n.\]

\end{thm}

\begin{thm}{Exercise}\label{ex:besikovitch=}
Show that if equality holds in \ref{thm:besikovitch+},
then $(M,g)$ is isometric to the rectangle $[0,d_1]\times\dots\times[0, d_n]$.
\end{thm}

\begin{thm}{Exercise}\label{ex:hexagon}
Suppose $g$ is a metric tensor on a regular hexagon $\textbf{\hexagon}$ such that $g$-distances between the opposite sides are at least $1$.
Is there a positive lower bound on $\area(\textbf{\hexagon},g)$?
\end{thm}

\begin{thm}{Exercise}\label{ex:cylinder}
Let $g$ be a Riemannian metric on the cylinder $\mathbb{S}^1\z\times [0,1]$.
Suppose that 
\begin{itemize}
\item 
$g$-distance between pairs of points on the opposite boundary circles $\mathbb{S}^1\times\{0\}$ and $\mathbb{S}^1\times\{1\}$ is at least 1, and 
\item
any curve $\gamma$ in $\mathbb{S}^1\times [0,1]$ that is homotopic to $\mathbb{S}^1\times\{0\}$ has $g$-length at least $1$.
\end{itemize}

\begin{subthm}{ex:cylinder:besicovitch}
Use Besicovitch inequality to show that
\[\area(\mathbb{S}^1\times [0,1],g)\ge \tfrac12.\]

\end{subthm}

\begin{subthm}{ex:cylinder:coarea}
Modify the proof of Besicovitch inequality using coarea inequality (\ref{cor:coarea}) to prove the optimal bound  
\[\area(\mathbb{S}^1\times [0,1],g)\ge 1.\]
 
\end{subthm}

\end{thm}

\begin{thm}{Exercise}\label{ex:gadograph}

\begin{subthm}{ex:gadograph-besikovitch}
Generalize \ref{thm:besikovitch+} to noncontinuous metric tensor $g$ described the following way:
there are two Riemannian metric tensors $g_1$ and $g_2$ on $M$ and a subset $V\subset M$ bounded by a Lipschitz hypersurface $\Sigma$ such that 
$g=g_1$ at the points in $V$ and $g=g_2$ otherwise.
\end{subthm}

\begin{subthm}{ex:gadograph-gadograph}
Use part \ref{SHORT.ex:gadograph-besikovitch} to prove the following: 
Let $V$ be a compact set in the $n$-dimensional Euclidean space $\EE^n$ bounded by a Lipschitz hypersurface $\Sigma$.
Suppose $g$ is a Riemannian metric on $V$ such that 
\[\dist{p}{q}{g}\ge\dist{p}{q}{\EE^n}\]
for any two points $p,q\in \Sigma$.
Show that
\[\vol(V,g)\ge \vol(V)_{\EE^n}.\]
\end{subthm}

\end{thm}

\begin{thm}{Exercise}\label{ex:involution-of-sphere}
Suppose that sphere with Riemannian metric $(\mathbb{S}^2,g)$ admits an involution $\iota$ such that $\dist{x}{\iota(x)}{g}\ge 1$.

Show that 
\[\area(\mathbb{S}^2,g)\ge \tfrac1{1000}.\]
Try to show that 
\[\area(\mathbb{S}^2,g)\ge \tfrac12,
\quad \area(\mathbb{S}^2,g)\ge 1,
\quad\text{or}\quad\area(\mathbb{S}^2,g)\ge \tfrac4\pi\]

\end{thm}

\begin{thm}{Advanced exercise}\label{ex:involution-of-3sphere}
Construct a metric tensor $g$ on $\mathbb{S}^3$ such that (1) $\vol(\mathbb{S}^3,g)$ arbitrarily small and (2) there is an involution $\iota\:\mathbb{S}^3\z\to \mathbb{S}^3$ such that $\dist{x}{\iota(x)}{g}\ge 1$ for any $x\in \mathbb{S}^3$.
\end{thm}

\begin{thm}{Exercise}\label{ex:GH-vol}
Let $g_1,g_2,\dots$, and $g_\infty$ be metrics on a fixed compact manifold $M$.
Suppose that $\distfun_{g_n}$ uniformly converges to $\distfun_{g_\infty}$ as functions on $M\times M\to\RR$.
Show that 
\[\liminf_{n\to\infty}\vol(M,g_n)\ge \vol(M,g_\infty).\]

Show that the inequality might be strict.
\end{thm}

\section{Systolic inequality}

Let $\spc{M}$ be a compact Riemannian space.
The \index{systole}\emph{systole} of $\spc{M}$ (briefly $\sys\spc{M}$) is defined to be the least length of a noncontractible closed curve in $\spc{M}$.

Let $\Lambda$ be a class of closed $n$-dimensional Riemannian spaces.
We say that a \index{systolic inequality}\emph{systolic inequality} holds for $\Lambda$ if there is a constant $c$ such that 
\[\sys\spc{M}\le c\cdot \sqrt[n]{\vol\spc{M}}\]
for any $\spc{M}\in \Lambda$.

\begin{thm}{Exercise}\label{ex:sysT2}
Use \ref{thm:besikovitch} or \ref{ex:cylinder} to show that a systolic inequality holds for any Riemannian metric on the 2-torus $\TT^2$.
\end{thm}

\begin{thm}{Exercise}\label{ex:sysRP2}
Use \ref{thm:besikovitch} to show that a systolic inequality holds for any Riemannian metric on  the real projective plane $\RP^2$.
\end{thm}

\begin{thm}{Exercise}\label{ex:sysSg}
Use \ref{thm:besikovitch+} to show that systolic inequality holds for any Riemannian metric on any closed surfaces of positive genus.
\end{thm}

\begin{thm}{Exercise}\label{ex:sysS2xS1}
Show that no systolic inequality holds for Riemannian metrics on $\mathbb{S}^2\times\mathbb{S}^1$.
\end{thm}

In the following lecture we will show that systolic inequality holds for many manifolds, in particular for torus of arbitrary dimension.

\section{Generalization}\label{sec:hausdorff-measure}

The following proposition follows immediately from the definitions of Hausdorff measure (Section \ref{sec:vol-haus}).

\begin{thm}{Proposition}\label{prop:bilip-measure}
Let $\spc{X}$ and $\spc{Y}$ be metric spaces, $A\subset \spc{X}$
and
 $f\: \spc{X}\to \spc{Y}$ be a $\Lip$-Lipschitz map. 
Then 
\[\haus_\alpha [f(A)]\le \Lip^\alpha\cdot\haus_\alpha\, A\]
for any $\alpha$.
\end{thm}

The following exercise provides a weak analog of the Besicovitch inequality that works for arbitrary metrics.

\begin{thm}{Exercise}\label{ex:besikovitch++}
Let $M$ be manifold with boundary and $\rho$ is a pseudometric on $M$.
Suppose $\partial M$ admits a degree 1 map to the surface of the $n$-dimensional cube $\square$;
denote by $d_1,\dots, d_n$ the $\rho$-distances between the inverse images of pairs of opposite faces of $\square$ in $M$.
Then 
\[\haus_n(M,\rho)\ge d_1\cdots d_n.\]
\end{thm}

Recall that in $n$-dimensional Riemannian spaces we have 
\[\tfrac{\omega_n}{2^n}\cdot\haus_n=\vol_n.\]
Note that $\tfrac{\omega_n}{2^n}<1$ if $n\ge 2$.
Therefore, the conclusion in \ref{ex:besikovitch++} is weaker than in \ref{thm:besikovitch+} (the assumptions are weaker as well).

One can redefine systolic inequality on $n$-dimensional manifolds using the Hausdorff measure $\haus_n$ instead of the volume.
It is straightforward to prove analogs of the exercises \ref{ex:sysT2}--\ref{ex:sysS2xS1} with this definition.

\begin{thm}{Exercise}\label{ex:2top-discs}
Suppose that two embedded $n$-disks $\Delta_1,\Delta_2$ in a metric space $\spc{X}$ have identical boundaries.
Assume that $\spc{X}$ is contractible and $\haus_{n+1}\spc{X}=0$.
Show that $\Delta_1=\Delta_2$.
\end{thm}

\section{Remarks}\label{sec:besicovitch-remarks}

The optimal constants in the systolic inequality are known only in the following three cases:
\begin{itemize}
\item For real projective plane $\RP^2$ the constant is $\sqrt{\pi/2}$ --- the equality holds for a quotient of a round sphere by isometric involution. The statement was proved by Pao Ming Pu \cite{pu}.\label{page:pu}
\item For torus $\TT^2$ the constant is $\sqrt{2}/\sqrt[4]{3}$ --- the equality holds for a flat torus obtained from a regular hexagon by identifying opposite sides; this is the so-called \index{Loewner's torus inequality}\emph{Loewner's torus inequality}.
\item For the Klein bottle $\RP^2\#\RP^2$  the constant is $\sqrt{\pi}/2^{3/4}$ --- the equality holds for a certain nonsmooth metric.
The statement was proved by Christophe Bavard \cite{bavard}.
\end{itemize}
The proofs of these results use the so-called {}\emph{uniformization theorem}   available in the 2-dimensional case only.
These proofs are beautiful, but they are too far from metric geometry.
A good survey on the subject is written by Christopher Croke and Mikhail Katz \cite{croke-katz}.

An analog of Exercise \ref{ex:GH-vol} with Hausdorff measure instead of volume does not hold for general metrics on a manifold.
In fact there is a nondecreasing sequence of metric tensors $g_n$ on $M$, such that (1) $\vol(M,g_n)<1$ for any $n$ and (2) $\distfun_{g_n}$ converges to a metric on $M$ with arbitrary large Hausdorff measure of any given dimension; such examples were constructed by Dmitri Burago, Sergei Ivanov, and David Shoenthal \cite{burago-ivanov-shoenthal}.

\chapter{Width and systole}

This lecture is based on a paper of Alexander Nabutovsky \cite{nabutovsky}.

\section{Partition of unity}

\begin{thm}{Proposition}\label{thm:part-unit}
 Let $\{V_i\}$ be a finite open covering of a compact metric space ${\spc{X}}$.
Then there are Lipschitz functions $\psi_i\:{\spc{X}}\z\to[0,1]$ such that (1) if $\psi_i(x)>0$, then $x\in V_i$ and (2) for any $x\in {\spc{X}}$ we have
$$\sum_i\psi_i(x)=1.$$

\end{thm}

A collection of functions $\{\psi_i\}$ that meets the conditions in \ref{thm:part-unit} is called 
a \index{partition of unity}\emph{partition of unity subordinate to the covering} $\{V_i\}$.

\parit{Proof.}
Denote by $\phi_i(x)$ the distance from $x$ to the complement of $V_i$;
that is,
$$\phi_i(x)=\distfun_{{\spc{X}}\backslash V_i}(x).$$
Note $\phi_i$ is $1$-Lipschitz
for any $i$
and $\phi_i(x)>0$ if and only if $x\in V_i$.
Since $\{V_i\}$ is a covering, we have that
$$\Phi(x)\df\sum_i\phi_i(x)>0\ \ \text{for any}\ \ x\in {\spc{X}}.$$
Since $\spc{X}$ is compact, $\Phi>\delta$ for some $\delta>0$.
It follows that $x\mapsto\tfrac1{\Phi(x)}$ is a bounded Lipschitz function. 

Set 
$$\psi_k(x)=\frac{\phi_k(x)}{\Phi(x)}.$$
Observe that by construction the functions $\psi_i$ meet the conditions in the proposition.
\qedsf

\section{Nerves}

Let $\{V_1,\dots,V_k\}$ be a finite open cover of a compact metric space $\spc{X}$.
Consider an abstract simplicial complex $\spc{N}$, with one vertex $v_i$ for each set $V_i$ such that a simplex with vertices $v_{i_1},\dots, v_{i_m}$ is included in $\spc{N}$ if 
the intersection $V_{i_1}\cap\dots\cap V_{i_m}$ is nonempty.
\begin{figure}[h!]
\vskip-0mm
\centering
\includegraphics{mppics/pic-1402}
\end{figure}
The obtained simplicial complex $\spc{N}$ is called the \index{nerve}\emph{nerve of the covering $\{V_i\}$}.
Evidently $\spc{N}$ is a finite simplicial complex ---
it is a subcomplex of a simplex with the vertices $\{v_1,\dots,v_k\}$.

Note that the nerve $\spc{N}$ has dimension at most $n$ if and only if the covering $\{V_1,\dots,V_k\}$ has \index{multiplicity of covering}\emph{multiplicity} at most $n+1$;
that is, any point $x\in\spc{X}$ belongs to
at most $n+1$ sets of the covering.

Suppose $\{\psi_i\}$ is  
a partition of unity subordinate to the covering $\{V_1,\dots,V_k\}$.
Choose a point $x\in {\spc{X}}$.
Note that the set
$$\{v_{i_1},\dots,v_{i_n}\}=\set{v_i}{\psi_i(x)>0}$$
form vertices of a simplex in $\spc{N}$.
Therefore 
$$\bm{\psi}\:x\mapsto \psi_1(x)\cdot v_1+\psi_2(x)\cdot v_2+\dots+\psi_k(x)\cdot v_n.$$
describes a Lipschitz map from ${\spc{X}}$ to the nerve $\spc{N}$ of $\{V_i\}$.
In other words, $\bm{\psi}$ maps a point $x$ to the point in $\spc{N}$ with \index{barycentric coordinates}\emph{barycentric coordinates} $(\psi_1(x),\dots,\psi_k(x))$.

Recall that the \index{star}\emph{star} of a vertex $v_i$ (briefly $\Star_{v_i}$) is defined as the union of the interiors of all simplicies that have $v_i$ as a vertex.
Recall that $\psi_i(x)>0$ implies $x\in V_i$.
Therefore we get the following:

\begin{thm}{Proposition}\label{prop:space->nerve}
Let $\spc{N}$ be a nerve of an open covering $\{V_1,\z\dots,V_k\}$ of a compact metric space $\spc{X}$.
Denote by $v_i$ the vertex of $\spc{N}$ that corresponds to $V_i$.

Then there is a Lipschitz map $\bm{\psi}\:\spc{X}\to\spc{N}$ such that $\bm{\psi}(V_i)\z\subset\Star_{v_i}$ for every $i$.
\end{thm}

\section{Width}

Suppose $A$ is a subset of a metric space $\spc{X}$.
The radius of $A$ (briefly $\rad A$) is defined as the least upper bound on the values $R>0$ such that $\oBall(x,R)\supset A$ for some $x\in \spc{X}$.

\begin{thm}{Definition}\label{def:width}
Let $\spc{X}$ be a metric space.
The \index{width}\emph{$n$-th width} of $\spc{X}$ (briefly $\width_n\spc{X}$) is the least upper bound on values $R>0$ such that $\spc{X}$ admits a finite open covering $\{V_i\}$ with multiplicity at most $n+1$ and $\rad V_i< R$ for each $i$.
\end{thm}

\parbf{Remarks.}

\begin{itemize} 
\item Observe that 
\[\width_0\spc{X}\ge\width_1\spc{X}\ge\width_2\spc{X}\ge\dots\]
for any compact metric space $\spc{X}$.
Moreover, if $\spc{X}$ is connected, then 
\[\width_0\spc{X}=\rad\spc{X}.\]

\item Usually width is defined using diameter instead of radius, but the results differ at most twice.
Namely, if $r$ is the $n$-th radius-width and $d$ --- the $n$-th diameter-width, then 
$r\le d\le 2\cdot r$.

\item Note that \index{Lebesgue covering dimension}\emph{Lebesgue covering dimension} of $\spc{X}$ can be defined as the least number $n$ such that $\width_n\spc{X}=0$.

\item Another closely related notion is the so-called \index{macroscopic dimension}\emph{macroscopic dimension on scale $R$};
it is defined as the  least number $n$ such that $\width_n\spc{X}<R$.
\end{itemize}

\begin{thm}{Exercise}\label{ex:macrodimension}
Suppose $\spc{X}$ is a compact metric space such that any closed curve $\gamma$ in $\spc{X}$ can be contracted in its $R$-neighborhood.
Show that macroscopic dimension of $\spc{X}$ on scale $100\cdot R$ is at most 1.

What about quasiconverse? That is, suppose a simply connected compact metric space $\spc{X}$ has macroscopic dimension at most 1 on scale $R$, is it true that any closed curve $\gamma$ in $\spc{X}$ can be contracted in its $100\cdot R$-neighborhood?
\end{thm}

The following exercise gives a good reason for the choice of term \index{width}\emph{width}; it also can be used as an alternative definition.

\begin{thm}{Exercise}\label{ex:width=suprad(inv)}
Suppose $\spc{X}$ is a compact metric space.
Show that $\width_n\spc{X}<R$ if and only if there is a finite $n$-dimensional simplicial complex $\spc{N}$ and a continuous map $\bm{\psi}\:\spc{X}\to \spc{N}$
such that 
\[\rad[\bm{\psi}^{-1}(s)]<R\]
for any $s\in \spc{N}$.
\end{thm}

\section{Riemannian polyhedrons}

A \index{Riemannian!polyhedron}\emph{Riemannian polyhedron} is defined as a finite simplicial complex with a metric tensor on each simplex such that the restriction of the metric tensor to a subsimplex coincides with the metric on the subsimplex.

The {}\emph{dimension} of a Riemannian polyhedron is defined as the largest dimension in its triangulation.
For Riemannian polyhedrons one can define length of curves and volume the same way as for Riemannian manifolds.

The obtained metric space will be called \emph{Riemannian polyhedron} as well.
A \index{triangulation}\emph{triangulation} of Riemannian polyhedron  will always be assumed to have the above property on the metric tensor.

Further we will apply the notion of width only to compact Riemannian polyhedrons.
If $\spc{P}$ is an $n$-dimensional Riemannian polyhedron, then 
we suppose that
\[\width\spc{P}\df\width_{n-1}\spc{P}.\]

Suppose that $\spc{P}$ is an $n$-dimensional Riemannian polyhedron;
in this case we will use short cut $\vol$ for $\vol_n$.
Let us define \index{volume profile}\emph{volume profile} of $\spc{P}$ as a function 
returning largest volume of $r$-ball in~$\spc{P}$;
that is, the volume profile of $\spc{P}$ is a function $\VolPro_{\spc{P}}\:\RR_+\to\RR_+$ defined by 
\[\VolPro_{\spc{P}}(r)\df \sup\set{\vol \oBall(p,r)}{p\in\spc{P}}.\]
Note that 
$r\mapsto \VolPro_{\spc{P}}(r)$ is nondecreasing  and
\[\VolPro_{\spc{P}}(r)\le\vol\spc{P}\]
for any $r$.
Moreover, if $\spc{P}$ is connected, then the equality $\VolPro_{\spc{P}}(r)=\vol\spc{P}$ holds
for $r\ge \rad \spc{P}$.

Note that if $\spc{P}$ is a connected 1-dimensional Riemannian polyhedron, then 
\[\width\spc{P}=\width_0\spc{P}=\rad\spc{P}.\]

\begin{thm}{Exercise}\label{ex:1D-case}
Let $\spc{P}$ be a 1-dimensional Riemannian polyhedron.
Suppose that $\VolPro_{\spc{P}}(R)<R$ for some $R>0$.
Show that 
\[\width \spc{P}<R.\]
Try to show that $c=\tfrac 12$ is the optimal constant for which the following inequality holds: 
\[\width \spc{P}<c\cdot R.\]
\end{thm}

\section{Volume profile bounds width}

\begin{thm}{Theorem}\label{thm:width<volpro}
Let $\spc{P}$ be an $n$-dimensional Riemannian polyhedron. 
If the inequality 
\[R> n\cdot \sqrt[n]{\VolPro_{\spc{P}}(R)}\]
holds for {}\emph{some} $R>0$, then 
\[\width\spc{P}\le  R.\]
\end{thm}

Since $\VolPro_{\spc{P}}(R)\le \vol\spc{P}$ for any $R>0$,
we get the following:

\begin{thm}{Corollary}\label{thm:width<vol}
For any $n$-dimensional Riemannian polyhedron $\spc{P}$, we have
\[\width\spc{P}\le n\cdot \sqrt[n]{\vol\spc{P}}.\]

\end{thm}

The proof of \ref{thm:width<volpro} will be given at the very end of this section,
after discussing {}\emph{separating polyhedrons}. 

Let us start three technical statements.
The first statement can be obtained by modifying a smoothing procedure for functions defined on Euclidean space. 

A function $f$ defined on a Riemannian polyhedron $\spc{P}$ is called \index{piecewise smooth}\emph{piecewise smooth} if there is a triangulation of $\spc{P}$ such that restriction of $f$ to every simplex is smooth.

\begin{thm}{Smoothing procedure}\label{smoothing-procedure}
Let $\spc{P}$ be a Riemannian polyhedron and $f\:\spc{P}\to \RR$ be a 1-Lipschitz function.
Then for any $\delta>0$ there is a piecewise smooth 1-Lipschitz function $\tilde f\:\spc{P}\to \RR$ such that 
\[|\tilde f(x)-f(x)|<\delta\]
for any $x\in  \spc{P}$.
\end{thm}

The following statement can be proved by applying the classical Sard's theorem to each simplex of a Riemannian polyhedron.

\begin{thm}{Sard's theorem}\label{sard}\index{Sard's theorem}
Let $\spc{P}$ be an $n$-dimensional Riemannian polyhedron and $f\:\spc{P}\to \RR$ be a piecewise smooth function.
Then for almost all values $a\in\RR$, the inverse image $f^{-1}\{a\}$  is a Riemannian polyhedron of dimension at most $n-1$ (we assume that $f^{-1}\{a\}$ is equipped with the induced length metric).
\end{thm}

The following statement can be proved by applying the coarea inequality (\ref{cor:coarea}) to the restriction of $f$ to each simplex of the polyhedron and summing up the results.

\begin{thm}{Coarea inequality}\index{coarea inequality}\label{poly-coarea}
Let $\spc{P}$ be an $n$-dimensional Riemannian polyhedron and $f\:\spc{P}\to \RR$ be a piecewise smooth 1-Lipschitz function.
Set $v\z=\vol_n (f^{-1}[r,R])$ and $a(t)=\vol_{n-1}(f^{-1}\{t\})$.
Then 
\[\int_r^Ra(t)\cdot dt\ge v .\]
In particular there is a subset of positive measure $S\subset [r,R]$ such that the inequality 
\[a(t)\ge \frac v{R-r}\]
holds for any $t\in S$.
\end{thm}

\section*{Separating subpolyhedrons}

\begin{thm}{Definition}
Let $\spc{P}$ be an $n$-dimensional Riemannian polyhedron.
An $(n-1)$-dimensional subpolyhedron $\spc{Q}\subset\spc{P}$ is called \index{separating subpolyhedron}\emph{$R$-separating} if for each connected component $U$ of the complement $\spc{P}\backslash \spc{Q}$ we have 
\[\rad U<R.\]

\end{thm}

\begin{thm}{Lemma}\label{lem:separating}
Let $\spc{P}$ be an $n$-dimensional Riemannian polyhedron.
Then given $R>0$ and $\eps>0$ there is a $R$-separating subpolyhedron $\spc{Q}\subset\spc{P}$ such that for any $r_0<r_1\le R$ we have
\[\VolPro_{\spc{Q}}(r_0)< \tfrac1{r_1-r_0}\cdot \VolPro_{\spc{P}}(r_1)+\eps.\]

\end{thm}

The proof reminds the proof of the following statement about minimal surfaces: 
\emph{if a point $p$ lies on an compact area-minimizing surface $\Sigma$ and $\partial\Sigma \cap \oBall(p,r)=\emptyset$, then
\[\area(\Sigma\cap \oBall(p,r))\le \tfrac12\cdot \area\mathbb{S}^2\cdot r^2.\]
}

\parit{Proof.}
Choose a small $\delta>0$.
Applying the smoothing procedure (\ref{smoothing-procedure}), we can exchange each distance function $\distfun_p$ on $\spc{P}$ by $\delta$-close piecewise smooth 1-Lipschitz function, which will be denoted by $\widetilde \distfun_p$.

By Sard's theorem (\ref{sard}), for almost all values $c\z\in(r_0\z+\delta, r_1-\delta)$, the level set
\[\tilde S_c(p)=\set{x\in \spc{P}}{\widetilde \distfun_p(x)=c}\]
is a Riemannian polyhedron of dimension at most $n-1$.
Since $\delta$ is small, the coarea inequality (\ref{poly-coarea}) implies that $c$ can be chosen so that in addition the following inequality holds:
\begin{align*}
\vol_{n-1}\tilde S_c(p)&\le \tfrac1{r_1-r_0-2\cdot\delta}\cdot\vol_n[\oBall(p,r_1)]<
\\
&<\tfrac1{r_1-r_0}\cdot \VolPro_{\spc{P}}(r_1)+\tfrac\eps2.
\end{align*}

Suppose $\spc{Q}$ is an $R$-separating subpolyhedron in $\spc{P}$ with almost minimal volume;
say its volume is at most $\tfrac\eps2$-far from the greatest lower bound.
Note that cutting from $\spc{Q}$ everything inside $\tilde S_c(p)$ and adding $\tilde S_c(p)$ produces a $R$-separating subpolyhedron, say $\spc{Q}'$.%
\footnote{If $\dim\tilde S_c(p)<n-1$, then it might happen that $\dim\spc{Q}'<n-1$; so, by the definition, $\spc{Q}'$ is not separating.
It can be fixed by adding a tiny $(n-1)$-dimensional piece to $\spc{Q}'$.}

Since $\spc{Q}$ has almost minimal volume, we have
\[\vol_{n-1}[\spc{Q}\cap \oBall(p,r_0)_{\spc{P}}]-\tfrac\eps2\le \vol_{n-1}S_c(p).\]
Therefore 
\[\vol_{n-1}[\spc{Q}\cap \oBall(p,r_0)_{\spc{P}}]\le\tfrac1{r_1-r_0}\cdot \VolPro_{\spc{P}}(r_1)+\eps.
\eqlbl{eq:volQ<ProP}\]
Recall that $\spc{Q}$ is equipped with the induced length metric;
therefore $\dist{p}{q}{\spc{Q}}\ge \dist{p}{q}{\spc{P}}$ for any $p,q\in \spc{Q}$;
in particular, 
\[\oBall(p,r_0)_{\spc{Q}}\subset \spc{Q}\cap \oBall(p,r_0)_{\spc{P}}\]
for any $p\in \spc{Q}$ and $r_0\ge 0$.
Hence, \ref{eq:volQ<ProP} implies the lemma.
\qeds

\begin{thm}{Lemma}\label{lem:separating-width}
Let $\spc{Q}$ be an $R$-separating subpolyhedron in an $n$-dimensional Riemannian polyhedron $\spc{P}$.
Then 
\[\width\spc{Q}\le R
\quad\Longrightarrow\quad
\width\spc{P}\le R.\]
\end{thm}

\parit{Proof.}
Choose an open covering $\{V_1,\dots,V_k\}$ of $\spc{Q}$ as in the definition of width (\ref{def:width});
that is, it has multiplicity at most $n$ and $\rad V_i<R$ for any $i$. 

Note that $\{V_1,\dots,V_k\}$ can be converted into an open covering of
a small neighbourhood of $\spc{Q}$ in $\spc{P}$ without increasing the multiplicity.
This can be done by setting 
\[V_i'=\bigcup_{x\in V_i}\oBall(x,r_x),\]
where $r_x\df\tfrac1{10}\cdot\inf\set{\dist{x}{y}{}}{y\in \spc{Q}\backslash V_i}$.

By adding to  $\{V_i'\}$ all the components of $\spc{P}\backslash \spc{Q}$,
we increase the multiplicity by at most 1 and obtain a covering of $\spc{P}$.
The statement follows since $\dim \spc{P}= \dim \spc{Q}\z+1$.
\qeds

\section*{Proof assembling}

\parit{Proof of \ref{thm:width<volpro}.}
We apply induction on the dimension $n=\dim\spc{P}$.
The base case $n=1$ is given in \ref{ex:1D-case}.

Suppose that the  $(n-1)$-dimensional case is proved.
Consider an $n$-dimensional Riemannian polyhedron $\spc{P}$ and suppose
\[n\cdot \sqrt[n]{\VolPro\spc{P}(R)}< R\]
for some $R>0$.
Let $\spc{Q}$ be an $R$-separating subpolyhedron in $\spc{P}$ provided by \ref{lem:separating} for a small $\eps>0$.

Applying  \ref{lem:separating} for $r=\tfrac{n-1}n\cdot R$ and $R$, we have that 
\begin{align*}
\VolPro_\spc{Q}(r) &< \frac 1{R-r}\cdot \VolPro_\spc{P}(R)+\eps<
\\
&<\frac {n}{R}\cdot\left(\frac{R}{n}\right)^n=
\\
&=\left(\frac{r}{n-1}\right)^{n-1};
\end{align*}
that is, $(n-1)\cdot \sqrt[n-1]{\VolPro\spc{Q}(r)}< r$.
Since $\dim\spc{Q}= n-1$, by the induction hypothesis, we get that
\[\width\spc{Q}\le r<R.\]
It remains to apply \ref{lem:separating-width}.
\qeds

\section{Width bounds systole}

Recall that a topological space $K$ is called \index{aspherical space}\emph{aspherical} if any continuous map $\mathbb{S}^k\to K$ for $k\ge 2$ is null-homotopic.

\begin{thm}{Theorem}\label{thm:sys<width}
Suppose $\spc{M}$ is a compact aspherical $n$-dimensional Riemannian manifold.
Then 
\[\sys\spc{M}\le 6 \cdot \width \spc{M}.\]
\end{thm}

\begin{thm}{Lemma}\label{lem:aspherical-homotopy}
Let $K$ be an aspherical space and $\spc{W}$ a connected CW-complex.
Denote by $\spc{W}^k$ the k-skeleton of $\spc{W}$.
Then any continuous map $f\:\spc{W}^2\to K$ can be extended to a continuous map $\bar f\:\spc{W}\to K$

Moreover, if $p\in \spc{W}$ is a 0-cell and $q\in K$.
Then a continuous maps of pairs $\phi_0,\phi_1\:(\spc{W},p)\to(K,q)$ are homotopic if and only if $\phi_0$ and $\phi_1$ induce the same homomorphism on fundamental groups $\pi_1(\spc{W},p)\to\pi_1(K,q)$.
\end{thm}

\parit{Proof.}
Since $K$ is aspherical, any continuous map $\partial\mathbb{D}^n\to K$ for $n\ge 3$
is hull-homotopic;
that is, it can be extended to a map $\mathbb{D}^n\:\to K$.

It makes it possible to extend $f$ to $\spc{W}^3$, $\spc{W}^4$, and so on.
Therefore $f$ can be extended to whole $\spc{W}$.

The only-if part of the second part of lemma is trivial;
it remains to show the if part.

Sine $\spc{W}$ is connected, we can assume that $p$ forms the only 0-cell in $\spc{W}$;
otherwise we can collapse a maximal subtree of the 1-skeleton in $\spc{W}$ to $p$.
Therefore, $\spc{W}^1$ is formed by loops that generate $\pi_1(\spc{W},p)$.

By assumption, the restrictions of $\phi_0$ and $\phi_1$ to $\spc{W}^1$ are homotopic.
In other words the homotopy $\Phi\:[0,1]\times \spc{W}$ is defined on the 2-skeleton of $[0,1]\times \spc{W}$.
It remains to apply the first part of the lemma to the product $[0,1]\times \spc{W}$.
\qeds

\begin{thm}{Lemma}\label{lem:sys-homotopy}
Suppose $\gamma_0,\gamma_1$ are two paths between points in a Riemannian space $\spc{M}$ such that $\dist{\gamma_0(t)}{\gamma_1(t)}{\spc{M}}<r$ for any $t\in[0,1]$.
Let $\alpha$ be a shortest path from $\gamma_0(0)$ to $\gamma_1(0)$ and $\beta$ be a shortest path from $\gamma_0(1)$ to $\gamma_1(1)$. 
If $2\cdot r<\sys\spc{M}$, then there is a homotopy $\gamma_t$ from
$\gamma_0$ to $\gamma_1$ such that $\alpha(t)\equiv \gamma_t(0)$ and $\beta(t)\equiv \gamma_t(1)$.
\end{thm}

\parit{Proof.}
Set $s=\sys\spc{M}$; 
since $2\cdot r<s$, we have that $\eps=\tfrac1{10}(s-2\cdot r)>0$.

\begin{wrapfigure}{o}{34mm}
\vskip-0mm
\centering
\includegraphics{mppics/pic-1405}
\end{wrapfigure}

Note that we can assume that $\gamma_0$ and $\gamma_1$ are rectifiable;
if not we can homotopy each into a broken geodesic line kipping the assumptions true. 

Choose a fine partition $0\z=t_0\z<t_1\z<\z\dots\z<t_n=1$.
Consider a sequence of shortest paths $\alpha_i$ from $\gamma_0(t_i)$ to $\gamma_1(t_i)$.
We can assume that $\alpha_0=\alpha$, $\alpha_n=\beta$, and each arc $\gamma_j|_{[t_{i-1},t_i]}$ has length smaller than $\eps$.
Therefore, every quadrilateral formed by concatenation  of $\alpha_{i-1}$, $\gamma_1|_{[t_{i-1},t_i]}$, reversed $\alpha_i$, and reversed arc $\gamma_0|_{[t_{i-1},t_i]}$ has length smaller than $s$.
It follows that this curve is contractible.
Applying this observation for each quadrilateral, we get the statement.
\qeds

\parit{Proof of \ref{thm:sys<width}.}
Let $\spc{N}$ be the nerve of a covering $\{V_i\}$ of $\spc{M}$ and $\bm{\psi}\:\spc{M}\to\spc{N}$ be the map provided by \ref{prop:space->nerve}.
As usual, we denote by $v_i$ the vertex of $\spc{N}$ that corresponds to $V_i$.
Observe that $\dim\spc{N}<n$;
therefore, $\bm{\psi}$ kills the fundamental class of $\spc{M}$.

Let us construct a continuous map  $f\:\spc{N}\to  \spc{M}$ such that
$f\circ\bm{\psi}$ is homotopic to the identity map on $\spc{M}$.
Note that once $f$ is constructed, the theorem is proved.
Indeed, since $\bm{\psi}$ kills the fundamental class $[\spc{M}]$ of $\spc{M}$, so does $f\circ\bm{\psi}$.
Therefore, $[\spc{M}]=0$ --- a contradiction.

Set $R=\width \spc{M}$ and $s=\sys\spc{M}$.
Assume we choose $\{V_i\}$ as in the definition of width (\ref{def:width}).
For each $i$ choose a point $p_i\in \spc{M}$ such that $V_i\subset \oBall(p_i,R)$.

Set $f(v_i)=p_i$ for each $i$.
It defines the map $f$ on the 0-skeleton $\spc{N}^0$ of the nerve $\spc{N}$.
Further, $f$ will be defined step by step on the skeletons $\spc{N}^1,\spc{N}^2, \dots$ of $\spc{N}$.

Let us map each edge $[v_iv_j]$ in $\spc{N}$ to a shortest path $[p_ip_j]$.
It defines $f$ on $\spc{N}^1$.
Note that image of each edge is shorter than $2\cdot R$.

Suppose $[v_iv_jv_k]$ is a triangle in $\spc{N}$.
Note that perimeter of the triangle $[p_ip_jp_k]$ can not exceed $6\cdot R$.
Since $6\cdot R<s$, the contour of $[p_ip_jp_k]$ is contractible.
Therefore, we can extend $f$ to each triangle of~$\spc{N}$.
It defines the map $f$ on $\spc{N}^2$.

Finally, since $\spc{M}$ is aspherical, by \ref{lem:aspherical-homotopy}, the map $f$ can be extended to $\spc{N}^3$, $\spc{N}^4$ and so on.

It remains to show that $f\circ\bm{\psi}$ is homotopic to the identity map.
Choose a CW structure on $\spc{M}$ with sufficiently small cells, so that each cell lies in one of $V_i$.
Note that $\bm{\psi}$ is homotopic to a map $\bm{\psi}_1$ that sends $\spc{M}^k$ to $\spc{N}^k$ for any $k$.
Moreover, we may assume that (1) if a 0-cell $x$ of $\spc{M}$ maps to a $v_i$, then $x\in V_i$ and (2) each 1-cell  of $\spc{M}$ maps to an edge or a vertex of $\spc{N}$.
Choose a 1-cell $e$ in $\spc{M}$; by the construction, $f\circ\bm{\psi}_1$ maps $e$ to a shortest path $[p_ip_j]$ and $e$ lies $\oBall(p_i,R)$.
Observe that $[p_ip_j]$ is shorter than $2\cdot R$.
It follows that the distance between points on $[p_ip_j]$ and $e$ can not exceed $3\cdot R$.
Choose a shortest path $\alpha_i$ from every 0 cell $x_i$  of $\spc{M}$ to $p_j=f\circ\bm{\psi}_1(x_i)$.
It defines a homotopy on $\spc{M}^0$.
Since $6\cdot R<s$, \ref{lem:sys-homotopy} implies that this homotopy can be extended to $\spc{M}^1$.
By \ref{lem:aspherical-homotopy}, it can be extended to whole $\spc{M}$.
\qeds

\begin{thm}{Exercise}\label{ex:sys<width}
Analyze the proof of \ref{thm:sys<width} and improve its inequality to 
 \[\sys\spc{M}\le 4 \cdot \width \spc{M}.\]
\end{thm}

\begin{thm}{Exercise}\label{ex:fillrad-inj}
Modify the proof of \ref{thm:sys<width} to prove the following:

Suppose that $\spc{M}$ is a closed $n$-dimensional Riemannian manifold with \index{injectivity radius}\emph{injectivity radius} at least $r$; that is, if $\dist{p}{q}{\spc{M}}<r$, then a shortest path $[pq]_{\spc{M}}$ is uniquely defined.
Show that
\[\width\spc{M}\ge \tfrac{r}{2\cdot(n+1)}.\]

Use \ref{thm:width<vol} to conclude that $\vol\spc{M}\ge \eps_n \cdot r^n$
for some $\eps_n>0$ that depends only on $n$.
\end{thm} 

The second statement in the exercise is a theorem of Marcel Berger~\cite{berger-n};
an inequality with optimal constant (with equality for round sphere) was obtained by Marcel Berger and Jerry Kazdan \cite{berger-kazdan}.

\section{Essential manifolds}

To generalize \ref{thm:sys<width} further, we need the following definition.

\begin{thm}{Definition}\label{def:essential}
A closed manifold $M$ is called \index{essential manifold}\emph{essential} if it admits a continuous map $\iota\:M\to K$ to an aspherical CW-complex $K$ such that $\iota$ sends the fundamental class of $M$ to a nonzero homology class in $K$.
\end{thm}

Note that any closed aspherical manifold is essential --- in this case one can take $\iota$ to be the identity map on $M$.

The real projective space $\RP^n$ provides an interesting example of an essential manifold which is not aspherical.
Indeed, the infinite dimensional projective space $\RP^\infty$ is aspherical and for the natural embedding $\RP^n\hookrightarrow\RP^\infty$ the image $\RP^n$ does not bound in $\RP^\infty$.
The following exercise provides more examples of that type:

\begin{thm}{Exercise}\label{ex:connected-sum-essential}
Show that the connected sum of an essential manifold with any closed manifold is essential.
\end{thm}

\begin{thm}{Exercise}\label{ex:product-essential}
Show that the product of two essential manifolds is essential.
\end{thm}

Assume that the manifold $M$ is essential and $\iota \:M\to K$ as in the definition.
Following the proof of \ref{thm:sys<width}, we can homotope the map 
$f\circ\bm{\psi}\:M\to M$ to the identity on the 2-skeleton of $M$;
further since $K$ is aspherical, we can homotope the composition
$\iota\z\circ f\circ\bm{\psi}$ to  $\iota$. 
Existence of this extension implies that $\iota$ kills the fundamental class of $M$ --- a contradiction.
So, taking \ref{ex:sys<width} into account, we proved the following generalization of \ref{thm:sys<width}:

\begin{thm}{Theorem}\label{thm:sys<width++}
Suppose $\spc{M}$ is an essential Riemannian space.
Then 
\[\sys\spc{M}\le 4 \cdot \width \spc{M}.\]
\end{thm}

As a corollary from \ref{thm:sys<width++} and \ref{thm:width<vol} we get the so-called \index{systolic inequality!Gromov's systolic inequality}\emph{Gromov's systolic inequality}:

\begin{thm}{Theorem}\label{thm:sys+}
Suppose $\spc{M}$ is an essential $n$-dimensional Riemannian space.
Then 
\[\sys\spc{M}\le 4 \cdot n\cdot \sqrt[n]{\vol\spc{M}}.\]
\end{thm}

\section{Remarks}

Theorem \ref{thm:sys+} was proved originally by Mikhael Gromov \cite{gromov-1983} with a worse constant.
The given proof is a result of a sequence of simplifications given by Larry Guth \cite{guth},
Panos Papasoglu \cite{papasoglu},
Alexander Nabutovsky and Roman Karasev \cite{nabutovsky}.

The calculations could be done better; namely we could get
\[\width\spc{P}\le c_n\cdot \sqrt[n]{\vol\spc{P}},\]
where
$c_n=\sqrt[n]{n!/2}= \tfrac ne+o(n)$ \cite{nabutovsky}.
As a result, we may get a stronger statement in \ref{thm:sys+}:
\[\sys\spc{M}\le 4 \cdot c_n\cdot \sqrt[n]{\vol\spc{M}}.\]

A wide open conjecture says that for any $n$-dimensional essential manifold we have
\[\frac{\sys\spc{M}}{\sqrt[n]{\vol\spc{M}}}\le\frac{\sys\RP^n}{\sqrt[n]{\vol\RP^n}},\eqlbl{eq:RPn}\]
where we assume that the $n$-dimensional real projective space $\RP^n$ is equipped with a canonical metric.
In other words, the ratio in the right-hand side of \ref{eq:RPn} is the optimal constant in the Gromov's systolic inequality; this  ratio grows as $O(\sqrt n)$.
(The ratio for $n$-dimensional flat torus grows as $O(\sqrt n)$ as well.)

\appendix
\chapter{Semisolutions}
\parbf{\ref{ex:non-differentiable}.}
Choose a function $r\mapsto \alpha(r)$ such that $\alpha'(r)\cdot r\to 0$ and $\alpha(r)\to\infty$ as $r\to 0$.
Consider the reparametrization of the Euclidean plane given by $\iota\:(r,\theta)\mapsto (r,\theta+\alpha(r))$ in the polar coordinates.
Observe that $\iota$ is not differentiable at the origin, but the metric tensor $g$ induced by $\iota$  is continuous.

\medskip

For more on the subject read the paper of Eugenio Calabi and Philip Hartman \cite{calabi-hartman}. 

\parbf{\ref{ex:volume-preserving+short}};
\ref{SHORT.ex:volume-preserving+short:injective}.
Suppose $p=f(x)=f(y)$ and the points $x,y\in \spc{M}$ are distinct.
Since $f$ is short, we get for any $r>0$ the ball $\oBall(p,r)_{\spc{N}}$ contains the images of $\oBall(x,r)_{\spc{M}}$ and $\oBall(y,r)_{\spc{M}}$.
Since $f$ is volume-preserving, we get
\[
\vol\oBall(x,r)_{\spc{M}}
+
\vol\oBall(y,r)_{\spc{M}}
\le
\vol\oBall(p,r)_{\spc{N}}.
\eqlbl{vol+vol<vol}\]

By \ref{obs:lip-chart}, for any $\eps>0$ and all sufficiently small $r>0$ the volumes of the balls  $\oBall(x,r)_{\spc{M}}$, $\oBall(y,r)_{\spc{M}}$ and $\oBall(p,r)_{\spc{N}}$, lie in the range $\omega_n\cdot e^{\mp2\cdot n\cdot\eps}\cdot r^n$, where $\omega_n$ denotes the volume of the unit ball in the $n$-dimensional Euclidean space.
The latter contradicts \ref{vol+vol<vol} for appropriate choice of $\eps$ and $r$.

\parit{\ref{SHORT.ex:volume-preserving+short:bi}.}
Denote by $\sigma(r,a)$ the volume of union of two $r$-balls in the $n$-dimensional Euclidean space such that the distance between their centers is $a$.
Observe that the function $(a,r)\mapsto \sigma(r,a)$ is continuous and increasing in $a$ and $r$ for $a\le r$.
Further, note that
\[\sigma(\lambda\cdot r,\lambda\cdot a)=\lambda^n\cdot \sigma(r,a)\]
for any $\lambda>0$.

Choose a point $z\in \spc{M}$ and small $\eps>0$.
By \ref{obs:lip-chart} there is $R>0$ such that $\oBall(z,10\cdot R)$ admits a $e^{\mp\eps}$-bilipschitz map to the $n$-dimensional Euclidean space.

Choose $x,y\in \oBall(z, R)$.
The argument used in part \ref{SHORT.ex:volume-preserving+short:injective} implies that 
\[e^{-n\cdot\eps}\cdot \sigma(e^{-\eps}\cdot r, e^{-\eps}\cdot \dist{x}{y}{\spc{M}})
\le 
e^{n\cdot\eps}\cdot \sigma(e^{\eps}\cdot r, e^{\eps}\cdot \dist{f(x)}{f(y)}{\spc{N}}).
\eqlbl{eq:v(r,a)}\]
This inequality implies a lower bound on $\dist{f(x)}{f(y)}{\spc{N}}$ in terms of $\dist{x}{y}{\spc{M}}$.

Use the listed properties of the function $(a,r)\mapsto \sigma(r,a)$ to show that for any $c<1$ there is $\eps>0$ such that \ref{eq:v(r,a)} implies that $b>c\cdot a$ for all sufficiently small $a$.

Finally, since $\spc{M}$ and  $\spc{N}$ are length-metric spaces, part~\ref{SHORT.ex:volume-preserving+short:bi} implies that $f$ is locally distance preserving.
(An inclusion map from a nonconvex open subset to the plane gives an example of volume preserving short map that is not distance preserving.)

\medskip

A more general result is discussed by Paul Creutz and Elefterios Soultanis \cite{creutz-soultanis}.

\parbf{\ref{ex:compact-interior}.} Denote by $\spc{M}$ and $\spc{M}^\circ$ the space of $(M,g)$ and $(M^\circ,g)$;
further denote by $\bar{\spc{M}}^\circ$ the completion of $\spc{M}^\circ$.
Observe that the inclusion $M^\circ\hookrightarrow M$ induces a short onto map $\iota\:\bar{\spc{M}}^\circ\z\to\spc{M}$.

Recall that $M$ is bounded by hypersurface that is locally a graph.
Use it to show that any sufficiently short curve $\gamma$ in $(M,g)$ can be approximated by a curve in $\spc{M}^\circ$ with $g$-length arbitrary close to $\length_g\gamma$.
Conclude that $\iota$ is an isometry.

\parbf{\ref{ex:besikovitch=}.}
From the proof of Besicovitch inequality, one can see that the restriction of $\bm{f}$ to the interior of $\spc{M}$ is
(1) volume-preserving, and 
(2) its differential $d_p\bm{f}\:\T_p\to \T_{\bm{f}(p)}$ is an isometry for almost all $p$.

Since $\bm{f}$ is Lipschitz, (2) can be used to show that $\bm{f}$ is short.
It remains to apply \ref{ex:volume-preserving+short} and \ref{ex:compact-interior}.

\parbf{\ref{ex:hexagon}.}
Consider the hexagon with flat metric and curved sides shown on the diagram.
Observe that its area can be made arbitrarily small while keeping the distances from the opposite sides at least 1.

\begin{figure}[!ht]
\begin{minipage}{.48\textwidth}
\centering
\includegraphics{mppics/pic-27}
\end{minipage}\hfill
\begin{minipage}{.48\textwidth}
\centering
\includegraphics{mppics/pic-23}
\end{minipage}
\vskip-4mm
\end{figure}

\parbf{\ref{ex:cylinder};} \ref{SHORT.ex:cylinder:besicovitch}.
Let $\alpha$ be a shortest curve that runs between the boundary components of the cylinder.
Cut the cylinder along $\alpha$.
We get a square with Riemannian metric on it $(\square,g)$.

Two opposite sides of $\square$ correspond to the boundary components of the cylinder.
The other pair corresponds to the sides of the cut.
By assumption, the $g$-distance between the first pair of sides is at least 1.

Consider a shortest curve $\beta$ that connects this pair of sides;
let us keep the same notation for the projection of $\beta$ in the cylinder.

Note that a cyclic concatenation $\gamma$ of $\beta$ with an arc of $\alpha$ is homotopic to a boundary circle.
Therefore $\length_g\gamma\ge1$.
Since $\alpha$ is a shortest path, its arc cannot be longer than any curve connecting its ends; therefore 
\[\length_g\beta\ge \tfrac 12\cdot\length_g \gamma\ge \tfrac 12.\]
That is, the other pair of sides of $\square$ lies on $g$-distance at least $\tfrac12$ from each other.
By \ref{thm:besikovitch+}, $\area(\square,g)\ge \tfrac12$, hence the result.

\parit{\ref{SHORT.ex:cylinder:coarea}.}
Note that any curve in the cylinder that is bordant to a boundary component has length at least $1$.
Therefore if $0\le t\le  1$, then the level sets 
\[L_t=\set{x\in \mathbb{S}^1\times[0,1]}{\distfun_{\mathbb{S}^1\times\{0\}}(x)_{g}=t}\] have length at least $1$.
Applying the coarea inequality, we get that
\[\area(\mathbb{S}^1\times[0,1],g)\ge 1.\]

\parbf{\ref{ex:gadograph}}; \ref{SHORT.ex:gadograph-besikovitch}.
Argue the same way as in \ref{thm:besikovitch}, but observe in addition that $\vol \Sigma=\vol \bm{f}(\Sigma)=0$ and use it time to time.

\parit{\ref{SHORT.ex:gadograph-gadograph}.}
Without loss of generality, we may assume that $V$ lies in a unit cube~$\square$.
Consider a noncontinuous metric tensor $\bar g$ on $\square$ that coincides with $g$ inside $V$ and with the canonical flat metric tensor outside of~$V$.

Observe that the $\bar g$-distances between opposite faces of $\square$ are at least 1.
Indeed this is true for the Euclidean metric and the assumption $\dist{p}{q}{g}\ge\dist{p}{q}{\EE^d}$  guarantees that one cannot make a shortcut in~$V$.
Therefore, the $\bar g$-distances between every pair of opposite faces is at least as large as 1 which is the Euclidean distance.

Applying part \ref{SHORT.ex:gadograph-besikovitch}, we get that $\vol(\square,\bar g)\ge \vol\square$.
Whence the statement follows.

\parbf{\ref{ex:involution-of-sphere}.}
Let $x\in \mathbb{S}^2$ be a point that minimize the distance $|x-x'|_g$.
Consider a shortest path $\gamma$ from $x$ to $x'$.
We can assume that 
\[|x-x'|_g=\length \gamma=1.\]

Let $\gamma'$ be the antipodal arc to $\gamma$.
Note that $\gamma'$ intersects $\gamma$ only at the common endpoints $x$ and $x'$.
Indeed, if $p'=q$ for some $p,q\in\gamma$, then $|p-q|\ge 1$.
Since $\length \gamma=1$, the points $p$ and $q$ must be the ends of $\gamma$.

It follows that $\gamma$ together with $\gamma'$ forms a closed simple curve in $\mathbb{S}^2$;
it divides the sphere into two disks $D$ and $D'$.

Let us divide $\gamma$ into two equal arcs $\gamma_1$ and $\gamma_2$; each of length $\tfrac12$.
Suppose that $p,q\in\gamma_1$, then 
\begin{align*}
|p-q'|_g&\ge |q-q'|_g-|p-q|_g\ge
\\
&\ge 1-\tfrac12=\tfrac12.
\end{align*}
That is, the minimal distance from $\gamma_1$ to $\gamma_1'$ is at least~$\tfrac12$.
The same way we get that the minimal distance from $\gamma_2$ to $\gamma_2'$ is at least~$\tfrac12$.
By Besicovitch inequality, we get that 
\[\area(D,g)\ge \tfrac14\quad\text{and}\quad \area(D',g)\ge \tfrac14.\]
Therefore 
\[\area(\mathbb{S}^2,g)\ge\tfrac12.\]

\parit{A better estimate.}
Let us indicate how to improve the obtained bound to
\[\area(\mathbb{S}^2,g)\ge1.\]

Suppose $x$, $x'$, $\gamma$ and $\gamma'$ are as above.
Consider the function
\[f(z)=\min_t \{\,|\gamma'(t)-z|_g+t\,\}.\]
Observe that $f$ is 1-Lipschitz.

Show that two points $\gamma'(c)$ and $\gamma(1-c)$ lie on one connected component of the level set $L_c=\set{z\in\mathbb{S}^2}{f(z)=c}$;
in particular 
\[\length L_c\ge 2\cdot|\gamma'(c)-\gamma(1-c)|_g.\]
By the triangle inequality, we have that
\begin{align*}
|\gamma'(c)-\gamma(1-c)|_g&\ge 1-|\gamma(c)-\gamma(1-c)|_g=
\\
&=1-|1-2\cdot c|.
\end{align*}

The coarea inequality (\ref{cor:coarea})
\[\area(\mathbb{S}^2,g)\ge \int\limits_0^1\length L_c\cdot dc\]
finishes the proof.

The bound $\tfrac12$ was proved by Marcel Berger \cite{berger}. 
Christopher Croke conjectured that the optimal bound is $\tfrac4\pi$ and the round sphere is the only space that achieves this \cite[Conjecture 0.3 in][]{croke} --- if you solved the last part of the problem, then publish the result.

\begin{wrapfigure}{r}{20 mm}
\vskip-0mm
\centering
\includegraphics{mppics/pic-1305}
\end{wrapfigure}

\parbf{\ref{ex:involution-of-3sphere}.}
Given $\eps>0$, construct a disk $\Delta$ in the plane with 
\[\length\partial \Delta<10\ \ \text{and}\ \ \area \Delta<\eps\]
that admits an continuous involution $\iota$ such that 
\[|\iota(x)-x|\ge 1\]
for any $x\in\partial \Delta$.

An example of $\Delta$ can be guessed from the picture;
the involution $\iota$ makes a length preserving half turn of its boundary $\partial \Delta$.

Take the product $\Delta\times \Delta\subset \EE^4$;
it is homeomorphic to the 4-ball.
Note that 
$$\vol_3[\partial(\Delta\times \Delta)]=2\cdot\area \Delta\cdot\length \partial \Delta<20\cdot\eps.$$
The boundary $\partial(\Delta\times \Delta)$ is homeomorphic to $\mathbb{S}^3$
and the restriction of the involution $(x,y)\z\mapsto (\iota(x),\iota(y))$ has the needed property.

It remains to smooth $\partial(\Delta\times \Delta)$ a  bit.

\parit{Remark.} This example is given by Christopher Croke \cite{croke}.
Note that according to \ref{thm:sys+}, 
the involution $\iota$ cannot be made isometric.

\parbf{\ref{ex:GH-vol}.}
Note that if $(M,g_\infty)$ is $e^{\mp\eps}$-bilipschitz to a cube, then applying Besicovitch inequality, we get that 
\[\liminf_{n\to\infty} \vol (M,g_n)\ge e^{-n\cdot \eps}\cdot\vol (M,g_\infty).\]

By the Vitali covering theorem, given $\eps>0$, we can cover the whole volume of $(M,g_\infty)$ by $e^{\pm\eps}$-bilipschitz cubes.
Applying the above observation and summing up the results, we get that 
\[\liminf_{n\to\infty} \vol (M,g_n)\ge e^{-n\cdot \eps}\cdot\vol (M,g_\infty).\]
The statement follows since $\eps$ is an arbitrary positive number.

To solve the second part of the exercise, start with $g_\infty$ and construct $g_n$ by  adding many tiny bubbles.
The volume can be increased arbitrarily with an arbitrarily small change of metric.

\parit{Remark.}
A more general result was obtained by Sergei Ivanov~\cite{ivanov-1997}.
Note that the statement does not hold true for Gromov--Hausdorff convergence.
In fact any compact metric space $\spc{X}$ can be GH-approximated by a Riemannian surface with an arbitrarily small area.
To show the latter statement, approximate $\spc{X}$ by a finite graph $\Gamma$, embed $\Gamma$ isometrically to the Euclidean space, and pass to the surface of its neighborhood.

\parbf{\ref{ex:sysT2}.}
Set $s=\sys(\TT^2,g)$.

Cut $\TT^2$ along a shortest closed noncontractible curve $\gamma$.
We get a cylinder $(\mathbb{S}^1,g)$ with a Riemannian metric on it.

Applying the argument in \ref{ex:cylinder:besicovitch}, we get that the $g$-distance between the boundary components is at least $\tfrac s2$.
Then \ref{ex:cylinder:besicovitch} implies that the area of torus is at least $\tfrac{s^2}2$.

\parit{Remark.}
The optimal bound is $\tfrac{\sqrt{3}}{2}\cdot s^2$; see  \ref{sec:besicovitch-remarks}.

\begin{wrapfigure}{r}{44 mm}
\vskip-4mm
\centering
\includegraphics{mppics/pic-25}
\end{wrapfigure}

\parbf{\ref{ex:sysRP2}.}
Set $s\z=\sys (\RP^2,g)$.
Cut $(\RP^2,g)$ along a shortest noncontractible curve $\gamma$.
We obtain $(\DD^2,g)$ --- a disc with metric tensor which we still denote by $g$.

Divide $\gamma$ into two equal arcs $\alpha$ and $\beta$.
Denote by $A$ and $A'$ the two connected components of the inverse image of $\alpha$.
Similarly denote by $B$ and $B'$ the two connected components of the inverse image of $\beta$.

Let $\gamma_1$ be a path from $A$ to $A'$;
map it to $\RP^2$ and keep the same notation for it.
Note that $\gamma_1$ together with a subarc of $\alpha$ forms a closed noncontractible curve in $\RP^2$.
Since $\length\alpha=\tfrac s2$, we have that $\length\gamma_1\ge \tfrac s2$.
It follows that the distance between $A$ and $A'$ in $(\DD^2,g)$ is at least $\tfrac s2$.
The same way we show that the distance between $B$ and $B'$ in $(\DD^2,g)$ is at least $\tfrac s2$.

Note that $(\DD^2,g)$ can be parameterized by a square with sides $A$, $B$, $A'$ and $B'$ and apply \ref{thm:besikovitch} to show that 
\[\area(\RP^2,g)=\area(\DD^2,g)\ge \tfrac14\cdot s^2.\]

\parit{Remark.}
The optimal bound is $\tfrac2 \pi\cdot s^2$; see  \ref{sec:besicovitch-remarks}.
In fact any Riemannian metric on the disc with the boundary globally isometric to a unit circle with angle metric has the area at least as large as the unit hemisphere.
It is expected that the same inequality holds for any compact surface with connected boundary (not necessarily a disc);
this is the so-called \index{filling area conjecture}\emph{filling area conjecture} \cite[it is mentioned Mikhael Gromov in 5.5.$\mathrm{B}'(\mathrm{e}')$ of][]{gromov-1983}.

\parbf{\ref{ex:sysSg}.} Cut the surface along a shortest noncontractible curve $\gamma$. 
We might get a surface with one or two components of the boundary.
In these two cases repeat the arguments in \ref{ex:sysRP2} or \ref{ex:sysT2} using \ref{thm:besikovitch+} instead of \ref{thm:besikovitch}.

\parbf{\ref{ex:sysS2xS1}.} Consider the product of a small 2-sphere with the unit circle.

\parbf{\ref{ex:besikovitch++}.}
Apply the same construction as in the original Besicovitch inequality, assuming that the target rectangle
$[0,d_1]\times\dots\times [0,d_n]$ equipped with the metric induced by the $\ell^\infty$ norm;
apply \ref{prop:bilip-measure} where it is appropriate.

\parbf{\ref{ex:2top-discs}.} Suppose that $\Delta_1\ne\Delta_2$.
Consider the map $f\:\mathbb{S}^n\to \spc{X}$ such that the restriction to north and south hemispheres describe $\Delta_1$ and $\Delta_2$ respectively.
Show that if $\Delta_1\ne\Delta_2$, then $\mathbb{S}^n$ can be parameterized by the boundary of the unit cube $\square$ in such a way that for any pair $A$, $A'$ of opposite faces their images $f(A)$, $f(A')$ do not overlap.

Since $\spc{X}$ is contractible, the map $f$ can be extended to a map of the whole cube.
By \ref{ex:besikovitch++} 
\[\haus_{n+1}[f(\square)]>0,\]
a contradiction.

\parbf{\ref{ex:macrodimension}.}
The following claim resembles Besicovitch inequality;
it is key to the proof:
\begin{itemize}
 \item[$({*})$] Let $a$ be a positive real number.
 Assume that a closed curve $\gamma$ in a metric space $\spc{X}$ can be subdivided into 4 arcs $\alpha$, $\beta$, $\alpha'$, and $\beta'$ in such a way that 
 \begin{itemize}
 \item $|x-x'|>a$ for any $x\in\alpha$ and $x'\in \alpha'$
 and
 \item $|y-y'|>a$ for any $y\in\beta$ and $y'\in \beta'$.
 \end{itemize}
 Then $\gamma$ is not contractible in its $\tfrac a2$-neighborhood.
\end{itemize}

To prove $({*})$, consider two functions defined on $\spc{X}$ as follows:
\begin{align*}
w_1(x)&=\min \{\,a,\distfun_{\alpha}(x)\,\}
\\
w_2(x)&=\min \{\,a,\distfun_{\beta}(x)\,\}
\end{align*}
and the map $\bm{w}\:\spc{X}\to [0,a]\times[0,a]$, defined by
\[\bm{w}\:x\mapsto(w_1(x),w_2(x)).\]

Note that 
\begin{align*}
\bm{w}(\alpha)&=0\times [0,a],
&
\bm{w}(\beta)&=[0,a]\times 0,
\\
\bm{w}(\alpha')&=a\times [0,a],
&
\bm{w}(\beta')&=[0,a]\times a.
\end{align*} 
Therefore, the composition $\bm{w}\circ\gamma$ is a degree 1 map 
\[\mathbb{S}^1\to \partial([0,a]\times[0,a]).\] 
It follows that if $h\:\DD\to \spc{X}$ shrinks $\gamma$, then there is a point $z\in\DD$ such that 
$\bm{w}\circ h(z)=(\tfrac a2,\tfrac a2)$.
Therefore, $h(z)$ lies at distance at least $\tfrac a2$ from $\alpha$, $\beta$, $\alpha'$, $\beta'$
and therefore from $\gamma$.
It proves the claim.

\medskip

Coming back to the problem, let $\{W_i\}$ be an open covering of the real line with multiplicity $2$ and $\rad W_i<R$ for each $i$;
for example take the covering by the intervals $((i-\tfrac23)\cdot R,(i+\tfrac23)\cdot R)$.

Choose a point $p\in \spc{X}$.
Denote by $\{V_j\}$ the connected components of $\distfun_p^{-1}(W_i)$ for all $i$.
Note that $\{V_j\}$ is an open finite cover of $\spc{X}$ with multiplicity at most 2.
It remains to show that $\rad V_j<100\cdot R$ for each $j$.

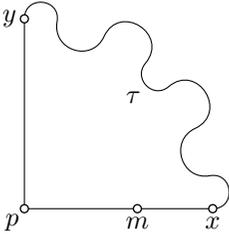
\begin{wrapfigure}{o}{31 mm}
\vskip-2mm
\centering
\includegraphics{mppics/pic-1310}
\end{wrapfigure}

Arguing by contradiction assume there is a pair of points  $x,y\in V_i$ 
such that $|x\z-y|_{\spc{X}}\z\ge 100\cdot R$.
Connect $x$ to $y$ with a curve $\tau$ in $V_j$.
Consider the closed curve $\sigma$ formed by $\tau$ and two shortest paths $[px]$, $[py]$.

Note that $|p-x|>40$.
Therefore, there is a point $m$ on $[px]$ such that $|m-x|=20$.

By the triangle inequality, the subdivision of $\sigma$ into the arcs $[pm]$, $[mx]$, $\tau$ and $[yp]$ satisfy the conditions of the claim $({*})$ for $a=10\cdot R$,
hence the statement.

\parit{The quasiconverse} does not hold.
As an example take a surface that looks like a long cylinder with closed ends;
\begin{figure}[h!]
\vskip0mm
\centering
\includegraphics{mppics/pic-1315}
\end{figure}
it is a smooth surface diffeomorphic to a sphere.
Assuming the cylinder is thin, it has macroscopic dimension 1 at a given scale.
However, a circle formed by a section of cylinder around its midpoint by a plane parallel to the base is a circle that cannot be contracted in its small neighborhood.

\parit{Source:} \cite[Appendix $1(\text{E}_{2})$]{gromov-1983}.

\parbf{\ref{ex:width=suprad(inv)}}; \textit{``only if'' part.}
Suppose $\width_n\spc{X}<R$.
Consider a covering $\{V_1,\dots,V_k\}$ of $\spc{X}$ guaranteed by the definition of width.
Let $\spc{N}$ be its nerve and $\bm{\psi}\:\spc{X}\to \spc{N}$ be the map provided by \ref{prop:space->nerve}.

Since the multiplicity of the covering is at most $n+1$, we have $\dim \spc{N}\le n$.

Note that if $x\in \spc{N}$ lies in a star of a vertex $v_i$,
then $\bm{\psi}^{-1}\{x\}\z\subset V_i$;
in particular, we have $\rad[\bm{\psi}^{-1}\{x\}]<R$.

\parit{``If'' part.}
Choose $x\in \spc{N}$.
Since the inverse image $\bm{\psi}^{-1}\{x\}$ is compact, $\bm{\psi}$ is continuous, and $\rad[\bm{\psi}^{-1}\{x\}]<R$,
there is a neighborhood $U\ni x$ such that the  $\rad[\bm{\psi}^{-1}(U)]<R$.

Since $\spc{X}$ is compact,  there is a finite cover $\{U_i\}$ of $\spc{N}$ such that $\bm{\psi}^{-1}(U_i)\subset\spc{X}$ has a radius smaller than $R$ for each $i$.
Since $\spc{N}$ has dimension $n$, we can inscribe%
\footnote{Recall that a covering $\{W_i\}$ is inscribed in the covering $\{U_i\}$ if for every $W_i$ is a subset of some $U_j$.} 
in $\{U_i\}$ a finite open cover $\{W_i\}$ with multiplicity at most $n+1$.
It remains to observe that $V_i=\bm{\psi}^{-1}(W_i)$ defines a finite open cover of $\spc{X}$ with  multiplicity at most $n+1$ and $\rad V_i<R$ for any $i$. 

\parbf{\ref{ex:1D-case}.}
Assume that $\spc{P}$ is connected.

Let us show that $\diam\spc{P}<R$.
If this is not the case, then there are points $p,q\in\spc{P}$ on distance $R$ from each other.
Let $\gamma$ be a shortest path from $p$ to $q$.
Clearly $\length\gamma\ge R$ and $\gamma$ lies in $\oBall(p,R)$ except for the endpoint $q$.
Therefore, $\length[\oBall(p,R)_{\spc{P}}]\ge R$.
Since $\VolPro_{\spc{P}}(R)\z\ge \length[\oBall(p,R)_{\spc{P}}]$,
the latter contradicts $\VolPro_{\spc{P}}(R)<R$.

In general case, we get that each connected component of $\spc{P}$ has a radius smaller than $R$.
Whence the width of $\spc{P}$ is smaller than $R$.

\parit{Second part.} Again, we can assume that $\spc{P}$ is connected.

The examples of line segment or a circle show that the constant $c=\tfrac12$ cannot be improved.
It remains to show that the inequality holds with $c=\tfrac12$.

Choose $p\in\spc{P}$ such that the value
\[\rho(p)=\max\set{\dist{p}{q}{\spc{P}}}{q\in\spc{P}}\]
is minimal.
Suppose $\rho(p)\ge\tfrac 12\cdot R$.
Observe that there is a point $x\z\in \spc{P}\backslash\{p\}$ that lies on any shortest path starting from $p$ and length $\ge\tfrac 12\cdot R$.
Otherwise for any $r\in(0,\tfrac 12\cdot R)$ there would be at least two points on distance $r$ from $p$;
by coarea inequality we get that the total length of $\spc{P}\cap \oBall(p,\tfrac 12\cdot R)$ is at least $R$ --- a contradiction.

Moving $p$ toward $x$ reduces $\rho(p)$ which contradicts the choice of~$p$.

\parbf{\ref{ex:sys<width}.}
The inequality $6\cdot R<s$ used twice:
\begin{itemize}
\item to shrink the triangle $[p_ip_jp_k]$ to a point;
\item to extend the constructed homotopy on $\spc{M}^0$ to $\spc{M}^1$.
\end{itemize}

The first issue can be resolved by passing to a barycentric subdivision of $\spc{N}^2$;
denote by $v_{ij}$ and $v_{ijk}$ the new vertices in the subdivision that correspond to edge $[v_iv_j]$ and triangle $[v_iv_jv_k]$ respectively.

Further for each vertex $v_{ij}$ choose a point $p_{ij}\in V_i\cap V_j$ and set $f(v_{ij})=p_{ij}$.
Similarly for each vertex $v_{ijk}$ choose a point $p_{ijk}\z\in V_i\cap V_j\cap V_k$ and set $f(v_{ijk})=p_{ijk}$.

Note that 
\[|p_i-p_{ij}|<R,\quad |p_i-p_{ijk}|<R,\quad\text{and}\quad |p_{ij}-p_{ijk}|<2\cdot R.\]
Therefore, perimeter of the triangle $[p_ip_{ij}p_{ijk}]$ in the subdivision is less that $4\cdot R$.
It resolves the first issue.

The second issue disappears if one estimates the distances a bit more carefully.
 
\parbf{\ref{ex:fillrad-inj}.}
Choose a fine covering of $\spc{M}$ with multiplicity at most $n$.
Choose $\bm{\psi}$ from $\spc{M}$ to the nerve $\spc{N}$ of the covering the same way as in the proof of \ref{thm:sys<width}.

It remains to construct $f\:\spc{N}\to\spc{M}$ and show that $f\circ\bm{\psi}$ is homotopic to the identity map.
To do this, apply the same strategy as in the proof of \ref{thm:sys<width} together with the so-called \index{geodesic cone construction}\emph{geodesic cone construction}
described below.

Let $\triangle$ be a simplex in a barycentric subdivision of $\spc{N}$.
Suppose that a map $f$ is defined on one facet $\triangle'$ of $\triangle$ to $\spc{M}$ and $\oBall(p,r)\supset f(\triangle')$.
Then one can extend $f$ to whole $\triangle$ such that the remaining vertex $v$ maps to $p$.
Namely connect each point $f(x)$ to $p$ by minimizing geodesic path $\gamma_x$ (by assumption it is uniquely defined) and set
\[f
\:
t\cdot x\z+(1-t)\cdot v
\mapsto
\gamma_x(t).\]

\parbf{\ref{ex:connected-sum-essential}.}
Suppose $M$ is an essential manifold and $N$ is an arbitrary closed manifold.
Observe that shrinking $N$ to a point produces a map $N\#M\to M$ of degree 1.
In particular, there is a map $f\:N\#M\to M$ that sends the fundamental class of $N\#M$ to the fundamental class of $M$.

Since $M$ is essential, there is an aspherical space $K$ and a map $\iota\:M\to K$ that sends the fundamental class of $M$ to a nonzero homology class in $K$.
From above, the composition $\iota\circ f\:N\#M\to K$ sends the fundamental class of $N\#M$ to the same homology class in~$K$.

\parbf{\ref{ex:product-essential}.}
Suppose $M_1$ and $M_2$ are essential.
Let $\iota_1\:M_1\to K_1$ and $\iota_2\:M_2\to K_2$ are the maps to aspherical spaces as in the definition (\ref{def:essential}).
Show that the map
$(\iota_1,\iota_2)\:M_1\times M_2\to K_1\times K_2$
meets the definition.

\parit{Remark.}
Choose a group $G$.
Note that there is an aspherical connected space CW-complex $K$ with fundamental group $G$.
The space $K$ is called an \index{K(G,1) space@$K(G,1)$ space}\emph{Eilenberg--MacLane space of type $K(G,1)$}, or briefly a $K(G,1)$ space.
Moreover it is not hard to check that
\begin{itemize}
\item $K$ is uniquely defined up to a weak homotopy equivalence;
\item if $\spc{W}$ is a connected finite CW-complex.
Then any homomorphism $\pi_1(\spc{W},w)\to\pi_1(K,k)$ is induced by a continuous map $\phi\:(\spc{W},w)\to(K,k)$.
Moreover, $\phi$ is uniquely defined up to homotopy equivalence.
\end{itemize}

\begin{itemize}
 \item Suppose that $M$ is a closed manifold, 
$K$ is a $K(\pi_1(M),1)$ space and a map $\iota\:M\to K$ induces an isomorphism of fundamental groups.
Then $M$ is essential if and only if $\iota$ sends the fundamental class of $M$ to a nonzero homology class of $K$.
\end{itemize}

The property described in the last statement is the original definition of essential manifold.
It can be used to prove a converse to the exercise;
namely \emph{the product of a nonessential closed manifold with any closed manifold is \emph{not} essential}.

%\include{tickets}
%%%%%%%%%%%%%%%%%%%%%%%%%%%%
{\small\sloppy
\input{metrics-on-manifolds.ind}

\printbibliography[heading=bibintoc]

@book {alexander-kapovitch-petrunin-2019,
    AUTHOR = {Alexander, S. and Kapovitch, V. and Petrunin, A.},
     TITLE = {An invitation to {A}lexandrov geometry: CAT(0) spaces},
      YEAR = {2019},
       URL = {https://doi.org/10.1007/978-3-030-05312-3},
}

@article {bavard,
    AUTHOR = {Bavard, C.},
     TITLE = {In\'{e}galit\'{e} isosystolique pour la bouteille de {K}lein},
   JOURNAL = {Math. Ann.},
  FJOURNAL = {Mathematische Annalen},
    VOLUME = {274},
      YEAR = {1986},
    NUMBER = {3},
     PAGES = {439--441},
      ISSN = {0025-5831},
   MRCLASS = {53C20},
  MRNUMBER = {842624},
MRREVIEWER = {R. Osserman},
       DOI = {10.1007/BF01457227},
       URL = {https://doi.org/10.1007/BF01457227},
}

@inproceedings {berger-kazdan,
    AUTHOR = {Berger, M. and Kazdan, J.},
     TITLE = {A {S}turm--{L}iouville inequality with applications to an
              isoperimetric inequality for volume in terms of injectivity
              radius, and to wiedersehen manifolds},
 BOOKTITLE = {General inequalities, 2 ({P}roc. {S}econd {I}nternat. {C}onf.,
              {O}berwolfach, 1978)},
     PAGES = {367--377},
 %PUBLISHER = {Birkh\"{a}user, Basel-Boston, Mass.},
      YEAR = {1980},
   MRCLASS = {53C20},
  MRNUMBER = {608261},
MRREVIEWER = {C. S. Houh},
       DOI = {10.1016/0377-0257(82)80029-3},
       URL = {https://doi.org/10.1016/0377-0257(82)80029-3},
}

@article {berger,
    AUTHOR = {Berger, M.},
     TITLE = {Volume et rayon d'injectivit\'{e} dans les vari\'{e}t\'{e}s riemanniennes de dimension {$3$}},
   JOURNAL = {Osaka Math. J.},
  FJOURNAL = {Osaka Mathematical Journal},
    VOLUME = {14},
      YEAR = {1977},
    NUMBER = {1},
     PAGES = {191--200},
      ISSN = {0388-0699},
   MRCLASS = {53C20},
  MRNUMBER = {467595},
MRREVIEWER = {R. L. Bishop},
       URL = {http://projecteuclid.org/euclid.ojm/1200770220},
}

@article {berger-n,
    AUTHOR = {Berger, M.},
     TITLE = {Volume et rayon d'injectivit\'{e} dans les vari\'{e}t\'{e}s riemanniennes},
   JOURNAL = {C. R. Acad. Sci. Paris S\'{e}r. A-B},
  FJOURNAL = {Comptes Rendus Hebdomadaires des S\'{e}ances de l'Acad\'{e}mie des
              Sciences. S\'{e}ries A et B},
    VOLUME = {284},
      YEAR = {1977},
    NUMBER = {19},
     PAGES = {A1221--A1224},
      ISSN = {0151-0509},
   MRCLASS = {53C20},
  MRNUMBER = {438251},
}

@Article{besicovitch,
    Author = {Besicovitch, A.},
    Title = {{On two problems of Loewner.}},
    FJournal = {{Journal of the London Mathematical Society}},
    Journal = {{J. Lond. Math. Soc.}},
    ISSN = {0024-6107; 1469-7750/e},
    Volume = {27},
    Pages = {141--144},
    Year = {1952},
    Publisher = {London Mathematical Society},
    Language = {English},
    DOI = {10.1112/jlms/s1-27.2.141},
    Zbl = {0046.05304}
}

@article{burago-ivanov-shoenthal,
    AUTHOR = {Burago, D. and Ivanov, S. and Shoenthal, D.},
     TITLE = {Two counterexamples in low-dimensional length geometry},
   JOURNAL = {Algebra i Analiz},
  FJOURNAL = {Rossi\u{\i}skaya Akademiya Nauk. Algebra i Analiz},
    VOLUME = {19},
      YEAR = {2007},
    NUMBER = {1},
     PAGES = {46--59},
      ISSN = {0234-0852},
   MRCLASS = {53C60},
  MRNUMBER = {2319509},
MRREVIEWER = {Andreas Bernig},
       DOI = {10.1090/S1061-0022-07-00984-3},
 }

@article {calabi-hartman,
    AUTHOR = {Calabi, E. and Hartman, P.},
     TITLE = {On the smoothness of isometries},
   JOURNAL = {Duke Math. J.},
  FJOURNAL = {Duke Mathematical Journal},
    VOLUME = {37},
      YEAR = {1970},
     PAGES = {741--750},
      ISSN = {0012-7094},
   MRCLASS = {53.74},
  MRNUMBER = {283727},
MRREVIEWER = {N. J. Hicks},
       URL = {http://projecteuclid.org/euclid.dmj/1077379303},
}

@article{croke,
    Author = {C. {Croke}},
    Title = {{Small volume on big $n$-spheres.}},
    FJournal = {{Proceedings of the American Mathematical Society}},
    Journal = {{Proc. Am. Math. Soc.}},
    ISSN = {0002-9939; 1088-6826/e},
    Volume = {136},
    Number = {2},
    Pages = {715--717},
    Year = {2008},
    Publisher = {American Mathematical Society (AMS), Providence, RI},
    Language = {English},
    DOI = {10.1090/S0002-9939-07-09079-X},
    MSC2010 = {53C20 52A40},
    Zbl = {1130.53023}
}

@incollection {croke-katz,
    AUTHOR = {Croke, C. and Katz, M.},
     TITLE = {Universal volume bounds in {R}iemannian manifolds},
 BOOKTITLE = {Surveys in differential geometry, {V}ol. {VIII} ({B}oston,
              {MA}, 2002)},
    SERIES = {Surv. Differ. Geom.},
    VOLUME = {8},
     PAGES = {109--137},
 PUBLISHER = {Int. Press, Somerville, MA},
      YEAR = {2003},
   MRCLASS = {53C23 (53C20 53C21 53C22)},
  MRNUMBER = {2039987},
MRREVIEWER = {Andrea Sambusetti},
       DOI = {10.4310/SDG.2003.v8.n1.a4},
       URL = {https://doi.org/10.4310/SDG.2003.v8.n1.a4},
}

@article {creutz-soultanis,
    AUTHOR = {Creutz, P. and Soultanis, E.},
     TITLE = {Maximal metric surfaces and the {S}obolev-to-{L}ipschitz
              property},
   JOURNAL = {Calc. Var. Partial Differential Equations},
  FJOURNAL = {Calculus of Variations and Partial Differential Equations},
    VOLUME = {59},
      YEAR = {2020},
    NUMBER = {5},
     PAGES = {Paper No. 177, 34},
      ISSN = {0944-2669},
   MRCLASS = {53C23 (30L10 49Q05 53B40)},
  MRNUMBER = {4153903},
       DOI = {10.1007/s00526-020-01843-0},
       URL = {https://doi.org/10.1007/s00526-020-01843-0},
}

@book {federer,
    AUTHOR = {Federer, H.},
     TITLE = {Geometric measure theory},
    SERIES = {Die Grundlehren der mathematischen Wissenschaften, Band 153},
 %PUBLISHER = {Springer-Verlag New York Inc., New York},
      YEAR = {1969},
    % PAGES = {xiv+676},
   MRCLASS = {28.80 (26.00)},
  MRNUMBER = {0257325},
MRREVIEWER = {J. E. Brothers},
}

@article {gromov-1983,
    AUTHOR = {Gromov, M.},
     TITLE = {Filling {R}iemannian manifolds},
   JOURNAL = {J. Differential Geom.},
  FJOURNAL = {Journal of Differential Geometry},
    VOLUME = {18},
      YEAR = {1983},
    NUMBER = {1},
     PAGES = {1--147},
      ISSN = {0022-040X},
   MRCLASS = {53C20 (53C21 57R99)},
  MRNUMBER = {697984},
MRREVIEWER = {Yu. Burago},
       URL = {http://projecteuclid.org/euclid.jdg/1214509283},
}

@article {guth,
    AUTHOR = {Guth, L.},
     TITLE = {Volumes of balls in {R}iemannian manifolds and {U}ryson width},
   JOURNAL = {J. Topol. Anal.},
  FJOURNAL = {Journal of Topology and Analysis},
    VOLUME = {9},
      YEAR = {2017},
    NUMBER = {2},
     PAGES = {195--219},
      ISSN = {1793-5253},
   MRCLASS = {53C23},
  MRNUMBER = {3622232},
MRREVIEWER = {St\'{e}phane Sabourau},
       DOI = {10.1142/S1793525317500029},
       URL = {https://doi.org/10.1142/S1793525317500029},
}

@article {ivanov-1997,
    AUTHOR = {Ivanov, S. V.},
     TITLE = {Gromov-{H}ausdorff convergence and volumes of manifolds},
   JOURNAL = {Algebra i Analiz},
  FJOURNAL = {Rossi\u{\i}skaya Akademiya Nauk. Algebra i Analiz},
    VOLUME = {9},
      YEAR = {1997},
    NUMBER = {5},
     PAGES = {65--83},
      ISSN = {0234-0852},
   MRCLASS = {53C23},
  MRNUMBER = {1604389},
MRREVIEWER = {Tadeusz Januszkiewicz},
}

@article {myers-steenrod,
    AUTHOR = {Myers, S. B. and Steenrod, N. E.},
     TITLE = {The group of isometries of a {R}iemannian manifold},
   JOURNAL = {Ann. of Math. (2)},
  FJOURNAL = {Annals of Mathematics. Second Series},
    VOLUME = {40},
      YEAR = {1939},
    NUMBER = {2},
     PAGES = {400--416},
      ISSN = {0003-486X},
   MRCLASS = {DML},
  MRNUMBER = {1503467},
       DOI = {10.2307/1968928},
       URL = {https://doi.org/10.2307/1968928},
}

@misc{nabutovsky,
    title={Linear bounds for constants in Gromov's systolic inequality and related results},
    author={A. Nabutovsky},
    year={2019},
    eprint={1909.12225},
    archivePrefix={arXiv},
    primaryClass={math.MG}
}

@misc{papasoglu,
    title={Uryson width and volume},
    author={P. Papasoglu},
    year={2019},
    eprint={1909.03738},
    archivePrefix={arXiv},
    primaryClass={math.DG}
}

@misc{petrunin2020pure,
    title={Pure metric geometry: introductory lectures},
    author={A. Petrunin},
    year={2020},
    eprint={2007.09846},
    archivePrefix={arXiv},
    primaryClass={math.MG}
}

@article {pu,
    AUTHOR = {Pu, P. M.},
     TITLE = {Some inequalities in certain nonorientable {R}iemannian
              manifolds},
   JOURNAL = {Pacific J. Math.},
  FJOURNAL = {Pacific Journal of Mathematics},
    VOLUME = {2},
      YEAR = {1952},
     PAGES = {55--71},
      ISSN = {0030-8730},
   MRCLASS = {53.0X},
  MRNUMBER = {48886},
MRREVIEWER = {A. Lichnerowicz},
       URL = {http://projecteuclid.org/euclid.pjm/1103051942},
}
\fussy
}

\end{document}